\documentclass[11pt]{article}

\usepackage[T1]{fontenc}
\usepackage{newtxtext}
\usepackage{mathtools,amssymb,amsthm,bm,mathrsfs}
\usepackage[a4paper,textwidth=16cm,textheight=23cm,centering]{geometry}
\usepackage{graphicx,tikz,float,hyperref,microtype}
\usepackage{enumitem}

\newenvironment{mathtable}[1]
  {\setlength{\arraycolsep}{4pt}\begin{array}{#1}}
  {\end{array}}

\DeclareMathSizes{10.5}{10.95}{7}{5} 
\hypersetup{colorlinks=true,linkcolor=black,urlcolor=blue!50!black}

\newtheorem{theorem}{Theorem}[section]
\newtheorem{lemma}[theorem]{Lemma}
\newtheorem{keylemma}[theorem]{Key Lemma}
\newtheorem{proposition}[theorem]{Proposition}

\theoremstyle{definition}
\newtheorem{definition}[theorem]{Definition}
\theoremstyle{remark}
\newtheorem{remark}[theorem]{Remark}
\newtheorem*{remark*}{Remark}

\numberwithin{equation}{section}
\mathtoolsset{showonlyrefs,showmanualtags=false}
\allowdisplaybreaks[1]

\title{\bfseries The Volume Conjecture for One-Edge 3-Manifolds with Geodesic Boundary}
\author{Yunpeng Meng \and Baojun Wu \and Yuxuan Yang}
\date{}

\begin{document}
\maketitle
\begin{abstract}
We prove the volume conjecture for compact orientable hyperbolic $3$-manifolds with totally geodesic boundary that admit a one-edge ideal triangulation; these manifolds form an infinite family. 
More precisely, we obtain an asymptotic formula for their Turaev--Viro invariant.
The proof uses new reflection symmetries of the quantum $6j$-symbol.
\end{abstract}

{\setlength{\parskip}{0pt}\tableofcontents}
\clearpage
\section{\texorpdfstring{Introduction}{Introduction}}

Volume conjectures relate the growth of quantum invariants at roots of
unity to hyperbolic volume. Kashaev \cite{Kashaev} formulated such a
conjecture for the link invariants constructed from the cyclic quantum
dilogarithm. For a hyperbolic link $K\subset S^3$, his conjecture
predicts that its $N$-th invariant $\langle K\rangle_N$ satisfies
$$
\lim\limits_{N\to\infty}\frac{2\pi}{N}\log|\langle K\rangle_N|
=\operatorname{Vol}(S^3\setminus K).
$$
Murakami and Murakami \cite{MurakamiMurakami} identified
$\langle K\rangle_N$ with the normalized $N$-colored Jones polynomial
evaluated at $e^{\frac{2\pi\mathbf i}{N}}$.

Chen and Yang \cite{CY} proposed analogous volume conjectures for
the Reshetikhin--Turaev and Turaev--Viro invariants of hyperbolic
$3$-manifolds. We consider their Turaev--Viro conjecture in the case
of manifolds with totally geodesic boundary.

\noindent\textbf{Chen--Yang volume conjecture}
\textup{(\cite[Conjecture~1.1]{CY}, totally geodesic boundary case).}
Let $M$ be a compact orientable hyperbolic $3$-manifold with nonempty
totally geodesic boundary. At $q=e^{\frac{2\pi\mathbf i}{r}}$, we have
$$
\lim\limits_{\substack{r\to\infty\\ r\ \mathrm{odd}}}
\frac{2\pi}{r}\log|\mathrm{TV}_r(M)|=\operatorname{Vol}(M).
$$

The doubling formula of Yang \cite[Theorem~1.2]{YangRelative}, together
with the asymptotic expansion conjecture of Wong and Yang
\cite[Conjecture~1.1]{WYLeading} applied to the closed double with
empty link, predicts the following refinement.

\noindent\textbf{Volume conjecture with Reidemeister torsion.}
Let $M$ be a compact orientable hyperbolic $3$-manifold with nonempty
totally geodesic boundary, and let $DM$ be its double along the boundary.
At $q=e^{\frac{2\pi\mathbf i}{r}}$, as $r\to\infty$ through odd integers,
$$
|\mathrm{TV}_r(M)|
=
\frac{2^{\dim_{\mathbb F_2}H_2(M;\mathbb F_2)}}{(4\pi)^{\chi(M)}}
\frac{r^{\frac{3}{2}\chi(M)}}{\sqrt{|\operatorname{Tor}(DM)|}}\,
e^{\frac{r}{2\pi}\operatorname{Vol}(M)}
\left(1+O\!\left(\frac1r\right)\right),
$$
where $\chi(M)$ is the Euler characteristic of $M$ and
$\operatorname{Tor}(DM)=\operatorname{Tor}(DM;\operatorname{Ad}\rho)$
is the Reidemeister torsion twisted by the adjoint action of the hyperbolic
holonomy $\rho\colon\pi_1(DM)\to\mathrm{PSL}(2,\mathbb C)$.

We consider compact orientable hyperbolic $3$-manifolds with totally
geodesic boundary that admit a one-edge ideal triangulation, a family
studied by Fujii~\cite{Fujii1990}, Kojima and Miyamoto~\cite{KojimaMiyamoto1991},
Frigerio, Martelli and Petronio~\cite{FMP,FMP2003Dehn,FMP2004Census},
Frigerio and Petronio~\cite{FP2004}, Chen and Yang~\cite{CY},
Jun Murakami~\cite{MurakamiKOS2020}, and Nigomedyanov~\cite{Nigomedyanov2026}.

A manifold $M_n$ admitting such a triangulation with $n\ge2$ tetrahedra
has an explicit hyperbolic structure: each tetrahedron is realized as
a regular truncated hyperideal tetrahedron with dihedral angle
$\frac{\pi}{3n}$~\cite[Remark~1.3]{FMP}. 
Since $\chi(M_n)=1-n$ and $H_2(M_n;\mathbb F_2)=0$ \cite[Lemma~1.1 and Theorem~1.2(4)]{FMP}, the homological factor in the conjecture is $1$ and the predicted power
of $r$ is $-\frac{3(n-1)}{2}$. We prove the refined conjecture for this family.

\begin{theorem}\label{thm:main}
Let $M_n$ be a compact orientable hyperbolic $3$-manifold with nonempty
totally geodesic boundary admitting an ideal triangulation with $n\ge2$
tetrahedra and a single edge. As $r\to\infty$ through odd integers,
at $q=e^{\frac{2\pi\mathbf i}{r}}$, we have
$$
\mathrm{TV}_r(M_n)
=(-1)^{\frac{(r+1)(n+1)}{2}}
\frac{(4\pi)^{n-1}}
{r^{\frac{3(n-1)}{2}}\sqrt{|\operatorname{Tor}(DM_n)|}}\,
e^{\frac{r}{2\pi}\operatorname{Vol}(M_n)}
\left(1+O\!\left(\frac1r\right)\right),
$$
where $\operatorname{Tor}(DM_n)$ is the Reidemeister torsion of the double $DM_n$
twisted by its adjoint holonomy.

In particular, the Chen--Yang volume conjecture holds for $M_n$:
$$
\lim\limits_{\substack{r\to\infty\\ r\ \mathrm{odd}}}
\frac{2\pi}{r}\log|\mathrm{TV}_r(M_n)|=\operatorname{Vol}(M_n).
$$
\end{theorem}

\paragraph{Outline of the proof.}
Section~2 describes the hyperbolic structure of $M_n$ and fixes the
normalization of the quantum $6j$-symbols and the Turaev--Viro invariant.

In Section~\ref{sec:alternative-sixj} we express two families of quantum
$6j$-symbol formulas, each with 18 members, as sums of the same
reflected kernel, with $k\in\mathbb Z$ in the first family and
$k\in b^{-2}+\mathbb Z$ in the second. These
reflection formulas constitute the main contribution of this paper.
Later in the proof, these symmetries play a central role in overcoming the
difficulty posed by a large Fourier term.

In Section~4 we use Poisson summation to express the Turaev--Viro
invariant in terms of Fourier integrals. Using the reflection formula of
Theorem~\ref{thm:hyperideal-finite-sums}, we first restrict the sum to the
symmetric region and then extend it from even $a$ to all integers $a$
(Key Lemmas~\ref{lem:symmetric-region} and~\ref{lem:reflection-extension}).
This avoids the central Fourier term whose potential function predicts a
growth rate greater than $\operatorname{Vol}(M_n)$. Meng, Wu, Yang and
Yang~\cite{MWWY} use a similar extension of the color sum in their study
of Reshetikhin--Turaev invariants.

In Section~5 we show that the imaginary part of the classical potential
at its geometric critical point is $\operatorname{Vol}(M_n)$ and relate
the Hessian determinant at this point to the Reidemeister torsion of
the double. We then compute the contribution of the four leading
Fourier terms and show that the remaining terms are exponentially
smaller.

Appendix~\ref{app:sixj-proofs} gives the square root conventions and proves
Proposition~\ref{prop:double-sine-sixj} and Theorem~\ref{thm:hyperideal-finite-sums}.

\noindent{\bf Acknowledgment.} 
We thank Tian Yang for helpful discussions.
B.W. is supported by National Key R\&D Program of China (No. 2023YFA 1010700) and NSF China Grant (No. 12526204).  Y.M. is partially supported by NSF Grant DMS-2505908.

\noindent{\bf AI Usage Disclosure.}
Our study of single-color reflections of quantum $6j$-symbols
was inspired by the $W_2(D_6)$ symmetry of the $b$-$6j$
symbols~\cite{TV}.
The authors used ChatGPT to  assist with the proof of Theorem~\ref{thm:hyperideal-finite-sums}.
The idea of using reflection symmetry to eliminate the problematic
Fourier terms originated with the authors.

 \section{\texorpdfstring{Preliminaries}{Preliminaries}}

\subsection{Hyperbolic \texorpdfstring{$3$}{3}-manifolds with totally geodesic boundary
admitting one-edge ideal triangulations}

For $n\ge2$, let $M_n$ be a connected, orientable compact
$3$-manifold with nonempty boundary, equipped with an ideal
triangulation $\mathcal T$ consisting of $n$ tetrahedra and a single
edge. Denote by $\mathcal M_n$ the set of homeomorphism classes of
such manifolds. This family was studied by Frigerio, Martelli and
Petronio~\cite{FMP}. The boundary of $M_n$ is connected and has
genus $n$~\cite[Lemma~1.1]{FMP}.

The number of homeomorphism classes in $\mathcal M_n$ has the following
growth type: there exist constants $0<c<C$ such that
\[
n^{cn}<\#\mathcal M_n<n^{Cn}
\]
for all sufficiently large $n$
\cite[Theorem~4.1]{FMP2003Dehn}.
In particular, these families contain infinitely many pairwise
nonhomeomorphic manifolds.
The first values are \cite[p.~428]{FMP2003Dehn}
\[
\#\mathcal M_2=8,\qquad
\#\mathcal M_3=74,\qquad
\#\mathcal M_4=2340,\qquad
\#\mathcal M_5=97568.
\]
Moreover, among the $5160$ compact orientable hyperbolic $3$-manifolds
with nonempty totally geodesic boundary admitting an ideal triangulation
with at most 4 tetrahedra, $2422$ admit a one-edge ideal triangulation
\cite[Section~2.5]{FMP2004Census}.

We can give $M_n$ a hyperbolic structure with totally geodesic boundary
as follows ~\cite[Remark~1.3]{FMP}. Realize each tetrahedron of $\mathcal T$ as a regular
truncated hyperideal tetrahedron with dihedral angle $\frac{\pi}{3n}$,
and glue the hexagonal faces by isometries according to the prescribed
face pairings. The total angle around the unique edge is
$6n\cdot\frac{\pi}{3n}=2\pi$, so the resulting metric has no cone
singularity. Since the truncation triangles meet the hexagonal faces
orthogonally, they glue together to form a totally geodesic boundary.

By rigidity for compact hyperbolic manifolds with totally geodesic
boundary~\cite{FP2004}, this hyperbolic structure is unique up to
isometry. In particular, its isometry class, and hence its volume,
is determined by the topology of $M_n$. If $\Delta_\theta$ denotes
the regular truncated hyperideal tetrahedron with dihedral angle
$\theta$, then
\begin{equation}\label{eq:one-edge-volume}
\operatorname{Vol}(M_n)
=n\operatorname{Vol}\!\left(\Delta_{\frac{\pi}{3n}}\right).
\end{equation}
Let $v_8=\operatorname{Vol}\!\left(\Delta_{0}\right)\approx 3.66386$ denote the volume of the regular ideal hyperbolic
octahedron. The volume estimate in
\cite[Theorem~1.1]{Nigomedyanov2026}, specialized to the case
without cusps, gives
\begin{equation}\label{eq:one-edge-volume-estimate}
nv_8-\frac{32\pi^2}{135n}
\le \operatorname{Vol}(M_n)
\le nv_8,
\qquad n\ge8.
\end{equation}
Moreover, these volumes are strictly increasing in $n$ for $n\ge2$
\cite[Corollary~1.3]{Nigomedyanov2026}.

we now describe the case $n=2 $. The 8 manifolds in $\mathcal M_2$, classified by
Fujii~\cite[Theorem~3.1]{Fujii1990}, are precisely the manifolds of
minimum volume among connected, orientable compact hyperbolic
$3$-manifolds with nonempty totally geodesic boundary
\cite{KojimaMiyamoto1991}. 
\begin{equation}\label{eq:mmin-volume}
\operatorname{Vol}(M_{2})
=2\operatorname{Vol}\!\left(\Delta_{\frac{\pi}{6}}\right)
\approx6.45199.
\end{equation}

We choose Fujii's example~$(a)$ as a representative.
The following description uses the face pairings in
\cite[Figures~3.7--3.8 and Table~3.10, column~(a)]{Fujii1990}.

Take two copies $T$ and $T'$ of $\Delta_{\frac{\pi}{6}}$,
with vertices before truncation labelled $0,1,2,3$ and
$0',1',2',3'$, respectively.
The corresponding truncation triangles carry the same numerical labels.
Label the hexagonal faces of $T$ opposite $0,1,2,3$
by $A,B,C,D$, respectively, and those of $T'$ opposite
$0',1',2',3'$ by $A,C,B,D$, respectively.
Label the 6 local edges of $T$ by
\[
a=01,\qquad b=02,\qquad c=03,\qquad
d=12,\qquad e=13,\qquad f=23,
\]
and use primed labels for the corresponding local edges of $T'$.

\begin{figure}[htbp]
\centering
\begin{minipage}[t]{0.50\textwidth}
\vspace{0pt}
\centering
\includegraphics[width=\linewidth]{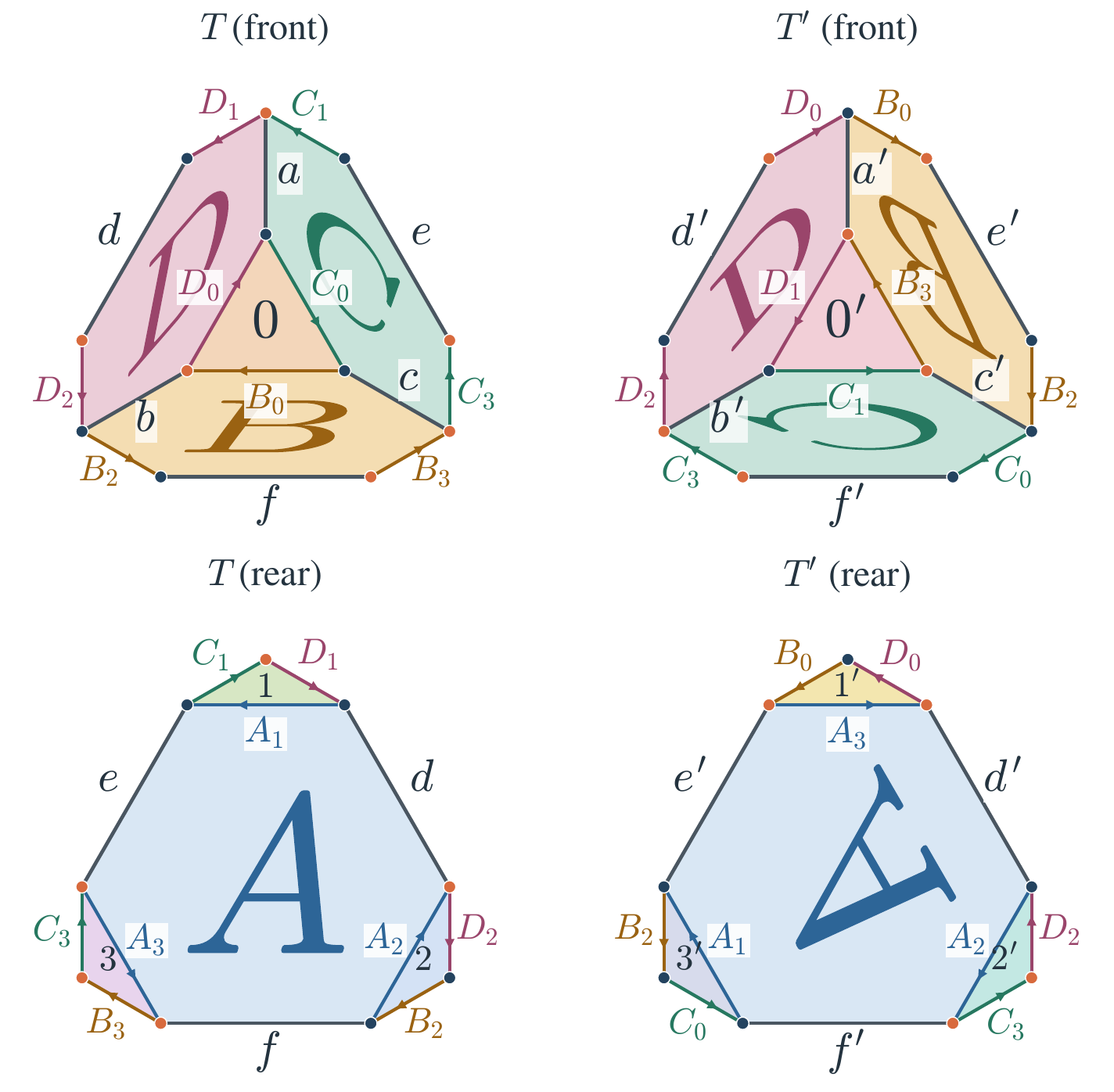}
\caption{Front and rear views of the tetrahedral face pairings.}
\label{f3}
\end{minipage}\hfill
\begin{minipage}[t]{0.49\textwidth}
\vspace{0pt}
\centering
\includegraphics[width=\linewidth]{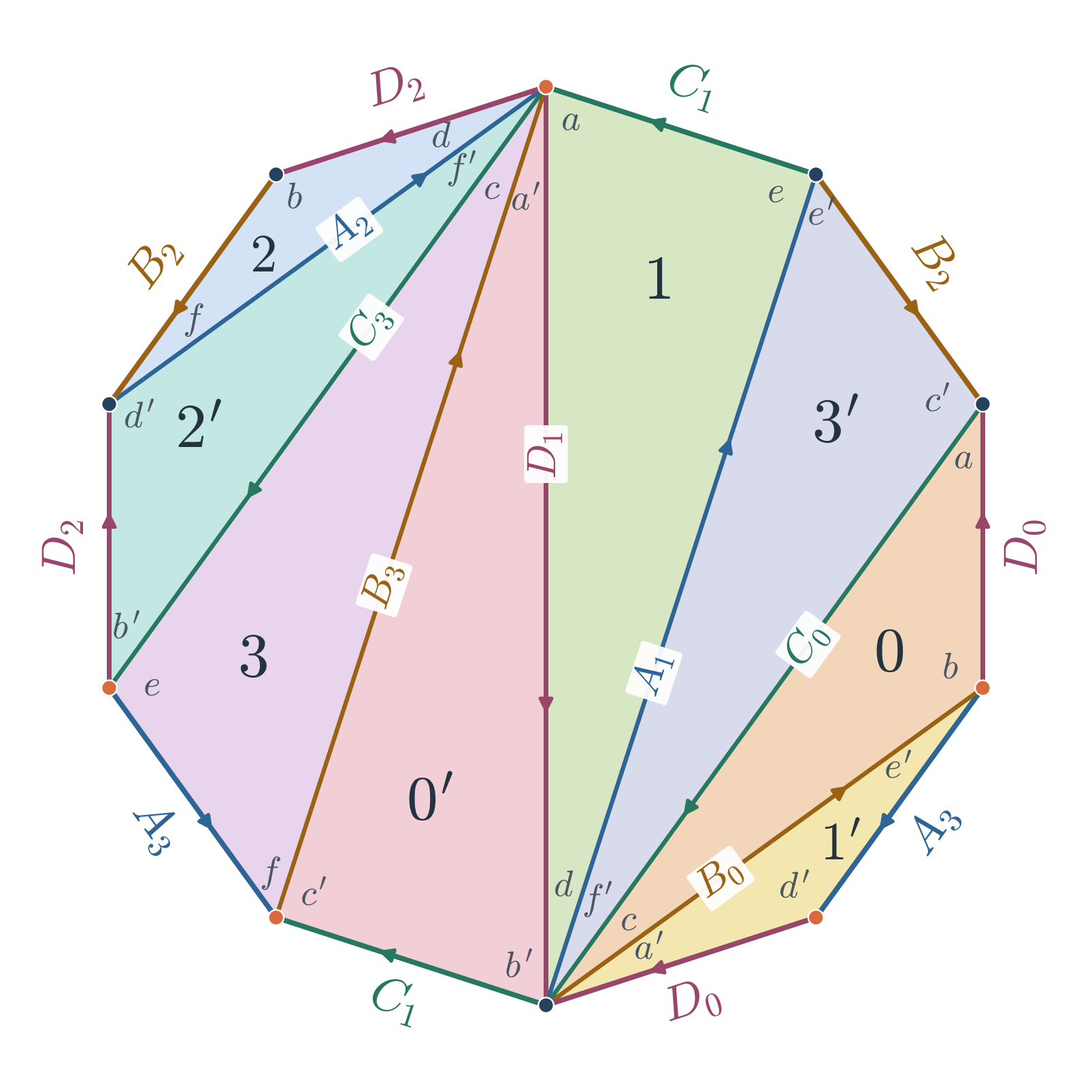}
\caption{A triangulation of the totally geodesic boundary.}
\label{fig:mmin-boundary}
\end{minipage}
\end{figure}

Identify faces carrying the same capital letter by the isometries
induced by the following ordered vertex correspondences:
\begin{equation}\label{eq:mmin-face-pairings}
\begin{mathtable}{c@{\qquad}c@{\;\longmapsto\;}c}
 A & (1,2,3) & (3',2',1')\\
 B & (0,2,3) & (1',3',0')\\
 C & (0,1,3) & (3',0',2')\\
 D & (0,1,2) & (1',0',2').
\end{mathtable}
\end{equation}
Each pairing reverses the induced orientations of the paired faces.

Label each side of a truncation triangle $i$ in $T$ by $X_i$,
where $X$ is the hexagonal face containing that side.
The corresponding side in $T'$ inherits the same label and arrow
under \eqref{eq:mmin-face-pairings}.

Corresponding elements in the two figures have matching labels and colors.
The letters at the corners of the boundary triangles indicate the
incident local edges.
In Figure~\ref{f3}, paired hexagonal faces carry the same capital letter.
The two copies of each letter are drawn to match under the face pairing
specified in \eqref{eq:mmin-face-pairings}.
These pairings identify all 12 local edges into a single edge
and all 8 vertices before truncation into one ideal vertex.
The 8 truncation triangles form its link and hence a connected
boundary surface.
The induced boundary triangulation has 2 distinct vertices,
corresponding to the two endpoints of the unique internal edge,
together with 12 edges and 8 triangles.
Its Euler characteristic is therefore $2-12+8=-2$, so its genus is~$2$.

\subsection{\texorpdfstring{Quantum $6j$-symbols}{Quantum 6j-symbols}}\label{subsec:quantum-sixj}

For $n\in\mathbb Z_{\ge0}$ the finite $q$-Pochhammer symbol is
\[
 (x;q)_n=\prod\limits_{i=0}^{n-1}(1-xq^i),\qquad (x;q)_0=1.
\]
For a positive integer $n$ and $q\ne0$, set $(x;q)_{-n}=\frac{1}{(xq^{-n};q)_n}$
wherever the denominator is nonzero.

Write $q=e^{\pi\mathbf i b^2}$.
For $\operatorname{Im}(b^2)>0$, the Pochhammer symbol of complex order $z$
is given by the convergent products
\[
 (x;q)_\infty=\prod\limits_{i=0}^{\infty}(1-xq^i),\qquad
 (x;q)_z=\frac{(x;q)_\infty}{(xq^z;q)_\infty}
\]
wherever the quotient is defined.
Define the quantum number and the generalized quantum factorial by
\begin{align}
 [z]&=q^{1-z}\frac{(q^{2z};q^2)_1}{1-q^2}
     =\frac{q^z-q^{-z}}{q-q^{-1}},\notag\\
 [z]!&=\frac{q^{-\frac{z(z-1)}{2}}}{(1-q^2)^z}
 \frac{(q^2;q^2)_z}
 {(e^{-\frac{2\pi\mathbf i}{b^2}};e^{-\frac{2\pi\mathbf i}{b^2}})_{-b^2z}} .
 \label{eq:pochhammer-factorial}
\end{align}
 We use the following conventions for complex powers:
$$
\begin{aligned}
q^z &:= e^{\pi i b^2 z},
&\qquad
(q^2)^z &:= e^{2\pi i b^2 z},\\
\left(e^{-2\pi i/b^2}\right)^z
&:= e^{-2\pi i z/b^2},
&\qquad
(1-q^2)^z
&:= \exp\!\bigl(z\log(1-q^2)\bigr).
\end{aligned}
$$
We use \eqref{eq:pochhammer-factorial} for
$\operatorname{Im}(b^2)>0$ and its value for real $b>0$ with
$b^2\notin\mathbb Z$, given by \eqref{eq:factorial-double-sine}.
This includes $b^2=\frac{2}{r}$, when
$q=e^{\frac{2\pi\mathbf i}{r}}$ for odd $r\ge3$. For an integer $n\ge0$,
$e^{2\pi\mathbf i n}=1$, so the two infinite products in the denominator
Pochhammer symbol cancel in their domain of convergence. Thus $[0]!=1$ and
\[
 [n]!=q^{-\frac{n(n-1)}{2}}\frac{(q^2;q^2)_n}{(1-q^2)^n}
      =\prod\limits_{i=1}^n[i].
\]
The finite product is a Laurent polynomial in $q$ and defines the usual
quantum factorial for every $q\ne0$, as well as for indeterminate $q$.

A triple of nonnegative integers is \emph{admissible} if it satisfies
the triangle inequalities and its sum is even.
A $6$-tuple $\mathbf{a}=(a_1,\ldots,a_6)$ is \emph{admissible} if the triples
$(a_1,a_2,a_3)$, $(a_1,a_5,a_6)$, $(a_2,a_4,a_6)$ and
$(a_3,a_4,a_5)$ are admissible.

For any $6$-tuple the face and quadrilateral parameters are
\begin{equation}\label{eq:T-parameters}
\begin{aligned}
T_1(\mathbf{a})&=\frac{a_1+a_2+a_3}{2},&
T_2(\mathbf{a})&=\frac{a_1+a_5+a_6}{2},\\
T_3(\mathbf{a})&=\frac{a_2+a_4+a_6}{2},&
T_4(\mathbf{a})&=\frac{a_3+a_4+a_5}{2},
\end{aligned}
\end{equation}
\begin{equation}\label{eq:Q-parameters}
\begin{aligned}
Q_1(\mathbf{a})&=\frac{a_1+a_2+a_4+a_5}{2},&
Q_2(\mathbf{a})&=\frac{a_1+a_3+a_4+a_6}{2},\\
Q_3(\mathbf{a})&=\frac{a_2+a_3+a_5+a_6}{2}.
\end{aligned}
\end{equation}
Set $Q_4=2b^{-2}-2$, independently of the colors.

For an integer $r\ge3$, an admissible triple $(a_i,a_j,a_k)$ is
\emph{$r$-admissible} if its entries lie in $\{0,\ldots,r-2\}$ and
$a_i+a_j+a_k\le 2(r-2)$. A $6$-tuple is \emph{$r$-admissible} if
each of its four face triples is $r$-admissible.

\begin{definition}\label{def:quantum-sixj}
For a $6$-tuple $\mathbf a$, define
\begin{equation}\label{eq:factorial-kernel}
 U_b(\mathbf a;k)=
 -\frac{1}{[Q_4+1]!}\,
 \frac{\sqrt{\displaystyle\prod\limits_{i,j=1}^4[Q_j-T_i]!}}
 {\displaystyle\prod\limits_{i=1}^4[k-T_i]!\prod\limits_{j=1}^4[Q_j-k]!}.
\end{equation}
Let $\mathbf a$ be admissible when $b^2\notin\mathbb Q$ and $r$-admissible
when $b^2=\frac{2}{r}$ with $r\ge3$ odd. Its \emph{quantum $6j$-symbol}
is the following sum over integer $k$:
\begin{equation}\label{eq:quantum-sixj}
 \left|\begin{matrix}a_1&a_2&a_3\\a_4&a_5&a_6\end{matrix}\right|
 =\sum\limits_{k=\max\{T_1,T_2,T_3,T_4\}}^{\min\{Q_1,Q_2,Q_3\}}
 U_b(\mathbf a;k).
\end{equation}
The numerator square root in \eqref{eq:factorial-kernel} is fixed by the
Racah convention in Appendix~\ref{app:root-signs}.
\end{definition}

\begin{remark*}
At $b^2=\frac{2}{r}$ we have $Q_4=r-2$. For $k\ge r-1$,
\eqref{eq:factorial-complement} gives
\[
 \frac{1}{[Q_4-k]!}=\frac{(-1)^{k+1}[k+1]!}{[r-1]!}=0,
\]
since $[k+1]!$ contains the factor $[r]=0$. Thus the terms with $k>Q_4$
in \eqref{eq:quantum-sixj} vanish.
\end{remark*}

The edge labeling is shown in Figure~\ref{f2}.


\subsection{Double sine function and dilogarithm function}
\label{subsec:special-functions}

To express the quantum $6j$-symbols in terms of the double sine
function, we recall its definition following
\cite[Sections~2.6--2.7]{MWWY} and \cite[Section~1.1]{MYComplex}.
Choose the square root $b$ with $\operatorname{Re}b>0$ and set $Q=b+b^{-1}$.
On the strip $0<\operatorname{Re}z<\operatorname{Re}Q$, define
$$
S_b(z)=
\exp\left(
\int_\Omega
\frac{\sinh\left(\left(\frac Q2-z\right)t\right)}
{4t\sinh\frac{bt}{2}
\sinh\frac{t}{2b}}\,dt
\right),
$$
where, for sufficiently small $\varepsilon>0$,
$$
\Omega=(-\infty,-\varepsilon]\cup
\{\varepsilon e^{\mathbf i\theta}:0\le\theta\le\pi\}\cup
[\varepsilon,\infty),
$$
with orientation from left to right. On this strip, we define
$$
\sqrt{S_b(z)}=
\exp\left(
\frac12\int_\Omega
\frac{\sinh\left(\left(\frac Q2-z\right)t\right)}
{4t\sinh\frac{bt}{2}
\sinh\frac{t}{2b}}\,dt
\right).
$$
The function $S_b$ extends meromorphically to $\mathbb C$, with
zeros at $Q+kb+\ell b^{-1}$ and poles at $-kb-\ell b^{-1}$,
where $k,\ell\in\mathbb{N}$.
\phantomsection\label{sb:functional-identities}By \cite[(2.1), (2.17)]{MYComplex}
and \cite[Lemma~3.2]{MWWY}, it satisfies
$$
\begin{cases}
S_b(z+b)&=2\sin(\pi bz)S_b(z)\\
S_b(z+b^{-1})&=2\sin\frac{\pi z}{b}\,S_b(z)
\end{cases},\quad S_b(z)S_b(Q-z)=1,\quad S_b(b)=b.
$$

The generalized quantum factorial of Section~\ref{subsec:quantum-sixj}
has the double sine expression
\begin{equation}\label{eq:factorial-double-sine}
 [z]!=(2\sin(\pi b^2))^{-z}\frac{S_b((z+1)b)}{S_b(b)} .
\end{equation}
With $Q_4=2b^{-2}-2$ and $n\in\mathbb Z_{\ge0}$, the complementary
factorial identities are
\begin{equation}\label{eq:factorial-complement}
 [z]!\left[\frac{Q_4}{2}-z\right]!=\left[\frac{Q_4}{2}\right]!,\qquad
 \frac{[Q_4+1]!}{[Q_4+1-n]!}=(-1)^n[n]!.
\end{equation}
Both formulas are proved in Appendix~\ref{app:special-function-identities}.

To estimate the summands as $r\to\infty$, we also recall the
dilogarithm function and the Lobachevsky function.
For $z\in\mathbb C\setminus (1,\infty)$, define
$$
\operatorname{Li}_2(z)=-\int_0^z\frac{\log(1-t)}{t}\,dt,
$$
using the principal branch of the logarithm. The Lobachevsky function is
$$
\Lambda(x)=-\int_0^x\log|2\sin t|\,dt,
\qquad x\in\mathbb R.
$$
It is odd and $\pi$-periodic, and
$$
|\Lambda(x)|\le\Lambda\left(\frac\pi6\right).
$$
For $0\le x\le \pi$, the two functions are related by
$$
\operatorname{Li}_2(e^{2\mathbf i x})
=\frac{\pi^2}{6}+x(x-\pi)+2\mathbf i\Lambda(x).
$$
For 
$x\in \mathbb{C}\backslash \{n\pi+ \mathbf{i} t |n\in \mathbb{Z},\ t<0\}$, set
$$
L(x)=x^2-\pi x+\frac{\pi^2}{6}
-\operatorname{Li}_2(e^{2\mathbf i x}).
$$
Then $L(x)=-2\mathbf i\Lambda(x)$ for $0\le x\le\pi$.
\begin{proposition}[{\cite[Proposition~2.15]{MWWY}}]
\label{prop:double-sine-asymptotics}
The following estimates are uniform for $0<b<1$, with $\log S_b$ given
by the defining integral.
\begin{enumerate}
\renewcommand{\labelenumi}{(\arabic{enumi})}
\item Define $\nu_b(x)$ by
$$
2\pi\mathbf i b^2\log S_b\left(\frac{x}{\pi b}+\frac b2\right)
-L(x)=\nu_b(x)b^4.
$$
For every $\varepsilon>0$, there is a constant $B_\varepsilon>0$
such that, for all $x\in\mathbb C$ with
$\varepsilon<\operatorname{Re}x<\pi-\varepsilon$,
$$
|\nu_b(x)|<B_\varepsilon.
$$
\item There is a constant $B>0$ such that, for all $x\in[0,\pi]$,
$$
\left|
2\pi\mathbf i b^2\log S_b\left(\frac{x}{\pi b}+\frac b2\right)
-L(x)
\right|<Bb^2.
$$
\end{enumerate}
\end{proposition}

We also use the following two estimates.

\begin{lemma}[{\cite[Lemma~2.1]{GMWY}}]
\label{lem:lobachevsky-continuity}
There is a constant $C>0$ such that, for all $x\in\mathbb R$
and $0<|a|<\frac12$,
\begin{equation}\label{eq:lobachevsky-continuity}
|\Lambda(x+a)-\Lambda(x)|\le C|a\log|a||.
\end{equation}
\end{lemma}

\begin{lemma}[{\cite[Lemma~2.6]{GMWY}}]
\label{lem:quantum-factorial-asymptotics}
Let $r\ge3$ be odd and $q=e^{\frac{2\pi\mathbf i}{r}}$. For every integer $0\le a\le r-1$,
\begin{equation}\label{eq:quantum-factorial-asymptotics}
\log|[a]!|
=-\frac{r}{2\pi}\Lambda\left(\frac{2\pi a}{r}\right)
-a\log\left(2\sin\frac{2\pi}{r}\right)+O(\log r),
\end{equation}
where the implicit constant is independent of $a$ and $r$.
\end{lemma}

For odd $r\ge3$, an $r$-admissible $6$-tuple of even colors is called
\emph{$r$-hyperideal} if, for each of the four admissible triples
and every permutation of its entries,
\[
0\le a_i+a_j-a_k< r-2,\qquad
r-2< a_i+a_j+a_k\le2(r-2).
\]

\begin{proposition}\label{prop:double-sine-sixj}
If $b^{-2}\notin\mathbb Z$, the kernel in \eqref{eq:factorial-kernel}
has the following double sine expression wherever its factors are defined.
\begin{equation}\label{eq:root-kernel}
 U_b(\mathbf a;k)=
 -\frac{b^2\sin(\pi b^2)}{\sin\!\left(\frac{\pi}{b^2}\right)}
 \frac{\sqrt{\displaystyle\prod\limits_{i,j=1}^4 S_b((Q_j-T_i+1)b)}}
 {\displaystyle\prod\limits_{i=1}^4 S_b((k-T_i+1)b)
  \prod\limits_{j=1}^4 S_b((Q_j-k+1)b)} .
\end{equation}
The square root is transported from \eqref{eq:factorial-kernel}.
\end{proposition}

Meng, Wu, Yang and Yang~\cite[Proposition~3.3]{MWWY} give the double sine sum
for the quantum $6j$-symbol in the special case of an $r$-hyperideal $6$-tuple
at the root of unity $q=e^{\frac{2\pi\mathbf i}{r}}$, with $r\ge3$ odd.
Appendix~\ref{app:double-sine-proof} proves Proposition~\ref{prop:double-sine-sixj}.
At $b^2=\frac{2}{r}$ its coefficient is
$\frac{2(-1)^{\frac{r+1}{2}}}{r}\sin\frac{2\pi}{r}$.

\subsection{Turaev--Viro invariants}\label{subsec:turaev-viro}

Fix an odd integer $r\ge3$ and take $q=e^{\frac{2\pi\mathbf i}{r}}$.

\begin{figure}[htbp]
    \centering
    \begin{minipage}[t]{0.4\textwidth}
        \centering
        \begin{minipage}[c][4.5cm][c]{\linewidth}
            \centering
            \includegraphics[width=\linewidth, height=4.5cm, keepaspectratio]{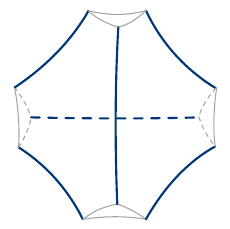}
        \end{minipage}
        \caption{Truncated tetrahedron}
        \label{f1}
    \end{minipage}
    \hfill
    \begin{minipage}[t]{0.59\textwidth}
        \centering
        \begin{minipage}[c][4.5cm][c]{\linewidth}
            \centering
\begin{tikzpicture}[
  x=0.0121875cm,y=-0.0121875cm,
  line width=.7pt,line cap=round,line join=round,
  font=\fontsize{10.5}{12.6}\selectfont
]
  \path[use as bounding box] (0,0) rectangle (640,420);
  \fill[white] (0,0) rectangle (640,420);
  \coordinate (T1) at (250,70);
  \coordinate (T2) at (90,340);
  \coordinate (T3) at (410,340);
  \coordinate (T4) at (550,105);
  \draw[dash pattern=on 2.8pt off 2.4pt] (T2)--(T4);
  \draw (T1)--(T2)--(T3)--(T4)--cycle;
  \draw[white,line width=2.5pt] (T1)--(T3);
  \draw (T1)--(T3);
  \foreach \p in {T1,T2,T3,T4}
    \fill (\p) circle[radius=1.17pt];
  \node at (246,32) {$v_1$};
  \node at (53,349) {$v_2$};
  \node at (435,381) {$v_3$};
  \node at (587,101) {$v_4$};
  \node at (140,205) {$a_1$};
  \node at (274,150) {$a_2$};
  \node at (402,52) {$a_3$};
  \node at (517,231) {$a_4$};
  \node at (401,150) {$a_5$};
  \node at (250,381) {$a_6$};
\end{tikzpicture}
        \end{minipage}
        \caption{The tetrahedral graph and its edge labels.}
        \label{f2}
    \end{minipage}
\end{figure}

Let $M$ be a compact $3$-manifold with non-empty boundary, and let
$\mathcal T$ be an ideal triangulation of $M$, that is, a finite
collection $\mathcal T=\{\Delta_1,\ldots,\Delta_{|\mathcal T|}\}$ of
truncated tetrahedra (see Figure \ref{f1}) with hexagonal faces identified in pairs by affine
homeomorphisms. We also let $E=\{e_1,\ldots,e_{|E|}\}$ be the set of
edges of $\mathcal T$.
A \emph{coloring} $\mathbf{a}$ of $(M,\mathcal T)$ assigns an integer
$a_i$ in between $0$ and $r-2$ to the edge $e_i$, and the coloring is
\emph{$r$-admissible} if the colors on the three edges of every face
form an $r$-admissible triple.

We use the normalization in \cite[Definition~1.1]{YangRelative}, with
the additional edge coloring there set equal to zero.

\begin{definition}
The $r$-th \emph{Turaev--Viro invariant} of $M$ is
\begin{equation}\label{eq:tv-invariant}
\mathrm{TV}_r(M)=
\sum\limits_{\mathbf{a}}
\left(\prod\limits_{i=1}^{|E|}(-1)^{a_i}[a_i+1]\right)
\prod\limits_{j=1}^{|\mathcal T|}
\left|\begin{matrix}
a_{j_1}&a_{j_2}&a_{j_3}\\
a_{j_4}&a_{j_5}&a_{j_6}
\end{matrix}\right|,
\end{equation}
where the sum is over all the $r$-admissible colorings
$\mathbf{a}=(a_1,\ldots,a_{|E|})$ of $(M,\mathcal T)$, and
$\{a_{j_1},\ldots,a_{j_6}\}$ are the colors of the edges of the
tetrahedron $\Delta_j$ assigned by $\mathbf{a}$ (see Figure \ref{f2}).
\end{definition}

The invariant in \eqref{eq:tv-invariant} is independent of the choice of
ideal triangulation \cite[Theorem~2.6]{CY}.

For a one-edged manifold $M_n$,
\begin{equation}\label{eq:one-edge-tv}
\mathrm{TV}_r(M_n)=
\sum\limits_{\substack{a\in2\mathbb{N}\\3a\le2(r-2)}}
[a+1]\left|\begin{matrix}a&a&a\\a&a&a\end{matrix}\right|^{\,n}.
\end{equation}
In particular, this invariant depends only on $n$ and $r$, and not
on the choice of $M_n\in\mathcal M_n$.

\section{Reflection symmetries of quantum \texorpdfstring{$6j$}{6j}-symbols}
\label{sec:alternative-sixj}

Recall that $Q_4=2b^{-2}-2$. At the root of unity
$q=e^{\frac{2\pi\mathbf i}{r}}$ with $r\ge3$ odd, we take
$b^2=\frac{2}{r}$, so $Q_4=r-2$.

Recall from \eqref{eq:quantum-sixj} that the sum over integer $k$ for
the original tuple $\mathbf a$ is
\[
 \left|\begin{matrix}a_1&a_2&a_3\\a_4&a_5&a_6\end{matrix}\right|
 =\sum\limits_{k=\max\{T_1,T_2,T_3,T_4\}}^{\min\{Q_1,Q_2,Q_3\}}
 U_b(\mathbf a;k).
\]

The single-color reflection $F_i$ replaces $a_i$ by $Q_4-a_i$ and fixes
the other five colors. For $i'=1,2,3$, the reflection $G_{i'}$ replaces each
of the four colors $a_s$ occurring in $Q_{i'}$ by $Q_{i'}(\mathbf a)-a_s$ and
fixes the other two.

Fix $i\in\{1,\ldots,6\}$. Denote the parameters omitting $a_i$ by
$T_{i_1},T_{i_2},Q_{i'_1}$ and those containing $a_i$ by
$T_{i_3},T_{i_4},Q_{i'_2},Q_{i'_3}$.

\begin{theorem}\label{thm:hyperideal-finite-sums}
Suppose either $b^2\notin\mathbb Q$ and $\mathbf a$ is admissible, or
$b^2=\frac{2}{r}$ for an odd integer $r\ge3$ and $\mathbf a$ is $r$-admissible.
For integer $k$ we have
\begin{equation}\label{eq:reflected-integer}
 \left|\begin{matrix}a_1&a_2&a_3\\a_4&a_5&a_6\end{matrix}\right|
 =\sum\limits_{k=\max\{T_{i_1}(\mathbf a),T_{i_2}(\mathbf a)\}}^{Q_{i'_1}(\mathbf a)}
 U_b(F_i(\mathbf a);k).
\end{equation}
Summing over $k\in b^{-2}+\mathbb Z$ also gives
\begin{equation}\label{eq:reflected-shifted}
 \left|\begin{matrix}a_1&a_2&a_3\\a_4&a_5&a_6\end{matrix}\right|
 =\sum\limits_{k=\max\{T_{i_3}(F_i(\mathbf a)),T_{i_4}(F_i(\mathbf a))\}}
       ^{\min\{Q_{i'_2}(F_i(\mathbf a)),Q_{i'_3}(F_i(\mathbf a))\}}
 U_b(F_i(\mathbf a);k).
\end{equation}
At $b^2=\frac{2}{r}$ the indices in the second formula are half-integers.

For either $i'\in\{i'_2,i'_3\}$, both identities also hold upon replacing
$\mathbf a$ by $G_{i'}(\mathbf a)$ throughout their right-hand sides.

Each of \eqref{eq:reflected-integer} and \eqref{eq:reflected-shifted}
thus gives 18 formulas, 6 from $F_i$ and 12 from $F_i\circ G_{i'}$.
\end{theorem}

For the proof, see Appendix~\ref{app:reflected-sums}.

\begin{remark}\label{rem:b-sixj-reflection}
For the $b$-$6j$ symbol, let $Q=b+b^{-1}$ and assume the entries are
$b$-admissible. The value of
\[
\left\{\begin{matrix}
\frac Q2\pm\frac{\theta_1}{2\pi b}&\frac Q2\pm\frac{\theta_2}{2\pi b}&\frac Q2\pm\frac{\theta_3}{2\pi b}\\[3pt]
\frac Q2\pm\frac{\theta_4}{2\pi b}&\frac Q2\pm\frac{\theta_5}{2\pi b}&\frac Q2\pm\frac{\theta_6}{2\pi b}
\end{matrix}\right\}_b
\]
is independent of each of the 6 signs by the reflection symmetry
of \cite[Proposition~2.3]{MYComplex}. If $\theta_s=\pi-\pi b^2(a_s+1)$, then $F_s$ sends $\theta_s$ to
$-\theta_s$. The reflected finite sums give a $q$-$6j$ counterpart of
this symmetry. For the weak reflection symmetry of extended quantum
$6j$-symbols, see \cite[Proposition~3.7]{MWWY}.
\end{remark}

\section{A simplified formula for the Turaev--Viro invariant}

Throughout Sections~4 and~5 we take $b=\sqrt{\frac{2}{r}}$ with $r\ge3$ odd.
Recall that $M_n$ has an ideal triangulation with $n\ge2$ tetrahedra and a single edge.
Using \eqref{eq:reflected-integer} of Theorem~\ref{thm:hyperideal-finite-sums},
Key Lemma~\ref{lem:symmetric-region}
restricts the even-color sum to the region $D_S$, invariant under
$x\mapsto-x$, with an exponentially smaller error.
In Section~\ref{subsec:reflection-extension}, Key Lemma~\ref{lem:reflection-extension}
extends the sum to all integers $a$.
This avoids the difficulty with the sum over even $a$, for which the potential function
of the central Fourier term predicts a growth rate greater than $\operatorname{Vol}(M_n)$.

\subsection{Reduction of the summation region}

By \eqref{eq:one-edge-tv}, the Turaev--Viro invariant of $M_n$ is
\[
\mathrm{TV}_r(M_n)=
\sum\limits_{\substack{a\in2\mathbb N\\3a\le 2(r-2)}}
[a+1]\left|\begin{matrix}a&a&a\\a&a&a\end{matrix}\right|^{n}.
\]
For $\mathbf a=(a,a,a,a,a,a)$, take $i=1$ in \eqref{eq:reflected-integer}.
Here $T_3=T_4=\frac{3a}{2}$ and $Q_3=2a$.
Since the terms with $k>r-2$ vanish, (3.1) gives
\[
\left|\begin{matrix}a&a&a\\a&a&a\end{matrix}\right|
=\sum\limits_{k=\frac{3a}{2}}^{\min\{2a,r-2\}} U_b(F_1(\mathbf a);k)
=\sum\limits_{k_i=\frac a2}^{\min\{a,r-2-a\}} U_b(F_1(\mathbf a);a+k_i).
\]
The bounds $\frac a2\le k_i\le r-2-a$ already imply $3a\le 2(r-2)$. Expanding the $n$-th power gives
\[
\mathrm{TV}_r(M_n)=
\sum\limits_{\substack{(a,\mathbf k)\in 2\mathbb N\times\mathbb Z^n\\
  \frac a2\le k_i\le\min\{a,r-2-a\},\ 1\le i\le n}}
[a+1]\prod\limits_{i=1}^{n} U_b(F_1(\mathbf a);a+k_i).
\]
Set
\[
x=\frac{2\pi}{r}(a+1)-\pi,\qquad
y_i=\frac{2\pi}{r}(k_i+1),\qquad 1\le i\le n.
\]
The lower bound $k_i\ge\frac a2$ becomes $y_i\ge\frac{\pi+x}{2}+\frac{\pi}{r}$.
The upper bounds $k_i\le a$ and $k_i\le r-2-a$ become $y_i\le\pi+x$ and
$y_i\le\pi-x$, giving $y_i\le\pi-|x|$. This gives the full summation region
\[
D_r=\left\{(x,\mathbf y)\in\mathbb R^{n+1}:
\frac{\pi+x}{2}+\frac{\pi}{r}\le y_i\le\pi-|x|,\quad 1\le i\le n\right\}.
\]

\begin{figure}[H]
\centering
\includegraphics[width=.75\linewidth]{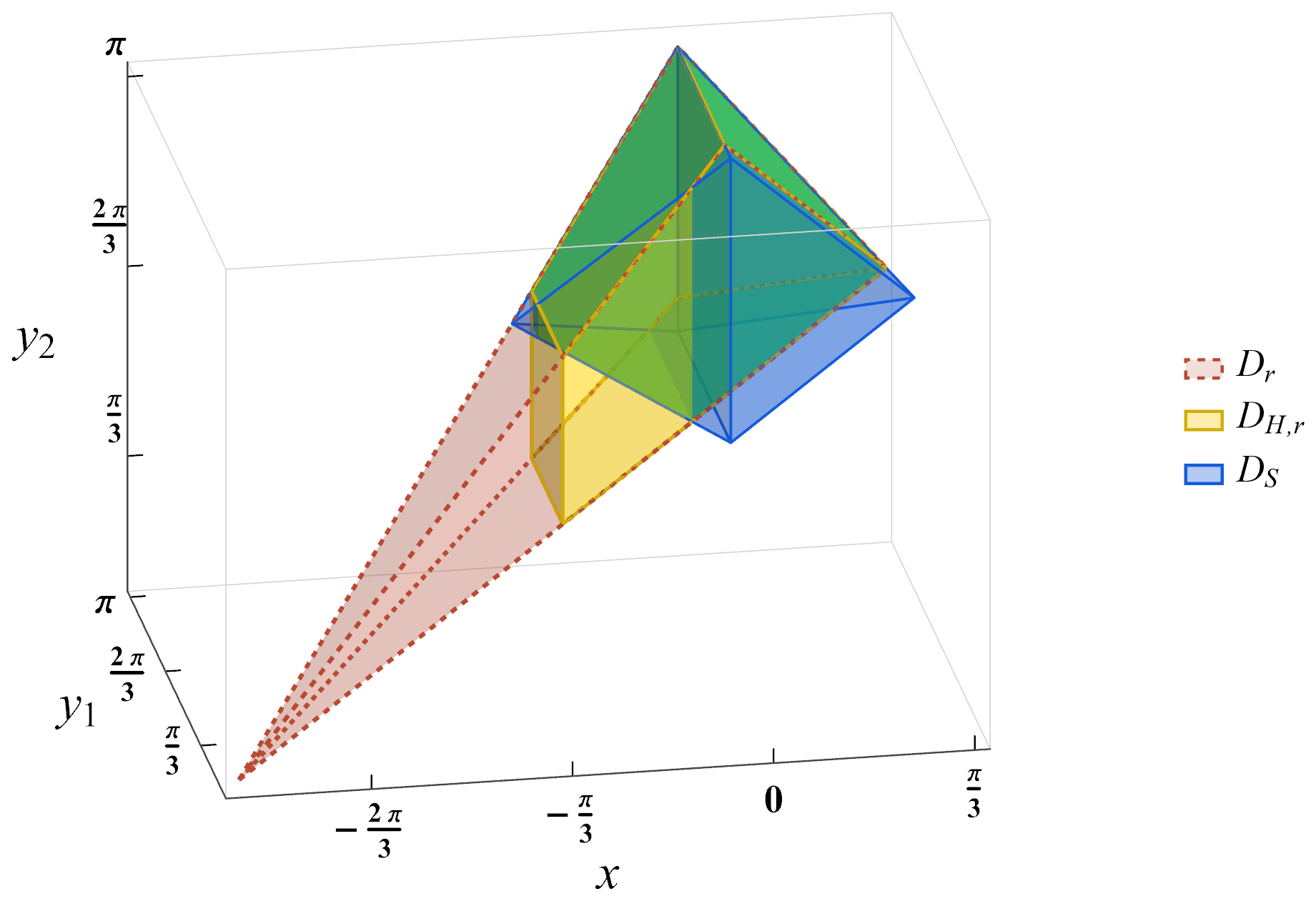}
\caption{The full region $D_r$, the hyperideal region $D_{H,r}$, and the symmetric region
$D_S$ for $n=2$ and $r=17$.}
\label{fig:summation-domains}
\end{figure}

Within $D_r$, the hyperideal condition reduces to $3a>r-2$, or equivalently
$x>-\frac{\pi}{3}+\frac{2\pi}{3r}$. This gives the hyperideal region
\[
D_{H,r}=\left\{(x,\mathbf y)\in D_r:
x>-\frac{\pi}{3}+\frac{2\pi}{3r}\right\}.
\]
We define the region $D_S$, symmetric under $x\mapsto -x$, by
\[
D_S=\left\{(x,\mathbf y)\in\mathbb R^{n+1}:
\frac{\pi+|x|}{2}<y_i<\pi-|x|,\quad 1\le i\le n\right\}.
\]
Our aim in this subsection is to reduce the summation region from $D_r$ to
$D_S$ up to an exponentially smaller error. Lemma~\ref{lem:hyperideal-colors}
first restricts the sum to $D_{H,r}$ by estimating the terms outside the
hyperideal range. Key Lemma~\ref{lem:symmetric-region} then bounds the
contribution from $D_{H,r}\setminus D_S$, leaving the sum over $D_S$.

\begin{lemma}\label{lem:one-edge-volume}
Let $M_n$ be as above.
\begin{enumerate}
\renewcommand{\labelenumi}{(\arabic{enumi})}
\item \cite[Theorem~1.2(1)]{FMP}
The hyperbolic volume of $M_n$ is
$$
\operatorname{Vol}(M_n)
=n\left[
8\Lambda\left(\frac{\pi}{4}\right)
-3\int_0^{\frac{\pi}{3n}}
\operatorname{arcosh}\left(\frac{\cos t}{2\cos t-1}\right)\,dt
\right].
$$
\item We have
$$
4n\Lambda\left(\frac{\pi}{6}\right)<\operatorname{Vol}(M_n).
$$
\end{enumerate}
\end{lemma}

\begin{proof}
Part (1) is the cited result.
For (2), since $n\ge2$ and the integrand is nonnegative,
$$
\frac{\operatorname{Vol}(M_n)}{n}
\ge8\Lambda\left(\frac{\pi}{4}\right)
-3\int_0^{\frac{\pi}{6}}
\operatorname{arcosh}\left(\frac{\cos t}{2\cos t-1}\right)\,dt.
$$
A numerical estimate shows that the right-hand side is greater
than $4\Lambda(\frac{\pi}{6})$, proving (2).
\end{proof}

When all 6 entries equal the same even integer $a$, the $r$-hyperideal condition is
$r-2<3a\le 2(r-2)$. The next lemma uses the volume bound above to
discard the terms outside this range and restrict the sum to $D_{H,r}$.

\begin{lemma}\label{lem:hyperideal-colors}
There exists $\varepsilon>0$, independent of $n$ and $r$, such that
for every $n\ge2$ and every $M_n$ as above,
$$
\mathrm{TV}_r(M_n)
=
\sum\limits_{\substack{a\in2\mathbb Z\\r-2<3a\le 2(r-2)}}
[a+1]
\left|\begin{matrix}a&a&a\\a&a&a\end{matrix}\right|^{n}
+
O\!\left(
e^{\frac{r}{2\pi}(\operatorname{Vol}(M_n)-\varepsilon)}
\right)
$$
as $r\to\infty$ through odd integers.
Here $O\!\left(e^{\frac{r}{2\pi}(\operatorname{Vol}(M_n)-\varepsilon)}\right)$
has an implied constant depending only on $n$.
\end{lemma}

\begin{proof}
Since $a$ is even and $r$ is odd, the terms omitted from
\eqref{eq:one-edge-tv} have $3a<r-2$. For these terms the Racah
formula \eqref{eq:racah-term} gives
$$
\left|\begin{matrix}a&a&a\\a&a&a\end{matrix}\right|
=
\sum\limits_{k=\frac{3a}{2}}^{2a}
\frac{(-1)^k\left[\frac{a}{2}\right]!^{6}[k+1]!}
{\left[\frac{3a}{2}+1\right]!^{2}
 \left[k-\frac{3a}{2}\right]!^{4}[2a-k]!^{3}}.
$$
Set
$$
\alpha=\frac{\pi a}{r},\qquad \beta=\frac{2\pi k}{r}-3\alpha.
$$
Then $0\le \beta\le \alpha<\frac{\pi}{3}$. By
Lemmas~\ref{lem:lobachevsky-continuity}
and~\ref{lem:quantum-factorial-asymptotics}, we obtain
$$
\left|
\log\left|
\frac{\left[\frac{a}{2}\right]!^{6}[k+1]!}
{\left[\frac{3a}{2}+1\right]!^{2}
 \left[k-\frac{3a}{2}\right]!^{4}[2a-k]!^{3}}
\right|
-\frac{r}{2\pi}F(\alpha,\beta)
\right|
\le C\log r,
$$
where $C>0$ is independent of $a,k,r$, and
$$
F(\alpha,\beta)
=-6\Lambda(\alpha)+2\Lambda(3\alpha)-\Lambda(3\alpha+\beta)
 +4\Lambda(\beta)+3\Lambda(\alpha-\beta).
$$

We claim that $F(\alpha,\beta)\le4\Lambda(\frac{\pi}{6})$.
If $0\le \alpha\le\frac{\pi}{6}$, then $\beta,\alpha-\beta\in[0,\alpha]$.
Since $\Lambda$ is increasing on $[0,\frac{\pi}{6}]$,
$$
4\Lambda(\beta)+3\Lambda(\alpha-\beta)\le7\Lambda(\alpha).
$$
Using $|\Lambda(x)|\le\Lambda(\frac{\pi}{6})$, we obtain
$$
F(\alpha,\beta)
\le\Lambda(\alpha)+2\Lambda(3\alpha)-\Lambda(3\alpha+\beta)
\le4\Lambda\left(\frac{\pi}{6}\right).
$$
If $\frac{\pi}{6}\le \alpha<\frac{\pi}{3}$, then $\Lambda(3\alpha)\le0$ and
$$
\Lambda(\alpha)\ge\Lambda\left(\frac{\pi}{3}\right)
=\frac{2}{3}\Lambda\left(\frac{\pi}{6}\right),
$$
where the equality follows from the duplication formula for $\Lambda$.
Hence
$$
F(\alpha,\beta)
\le-6\Lambda\left(\frac{\pi}{3}\right)
   +8\Lambda\left(\frac{\pi}{6}\right)
=4\Lambda\left(\frac{\pi}{6}\right).
$$

Thus, for each omitted color \(a\), the exponential growth rate of the
corresponding \(6j\)-symbol is at most \(4\Lambda(\pi/6)\). By Lemma~\ref{lem:one-edge-volume},
\[
\frac{\operatorname{Vol}(M_n)}n
\ge\frac{\operatorname{Vol}(M_2)}2=4\Lambda\!\left(\frac{\pi}{6}\right)+\varepsilon,
\qquad \varepsilon>0.
\]
For fixed $n$, the number of terms and the weights $|[a+1]|$ grow at most
polynomially in $r$, so the stated error bound follows.
\end{proof}

In the next key lemma, we use the new reflection symmetry of the quantum
$6j$-symbols in Theorem~\ref{thm:hyperideal-finite-sums}
to express the Turaev--Viro invariant in terms of a potential function that is
symmetric under $x\mapsto -x$.
We sum over $D_S$, which has the same symmetry, and show that the total contribution
from $D_{H,r}\setminus D_S$ is exponentially smaller than the leading term.

\begin{keylemma}\label{lem:symmetric-region}
Let $M_n$ be as above and let $a$ be an even integer.

For $\frac{\pi+|x|}{2}<y<\pi-|x|$, define the single-tetrahedron potential function
\[
\begin{aligned}
\tau_b(x,y)=&\pi\mathbf i b^2\Bigg[
 2\log S_b\!\left(\frac{y}{\pi b}\right)
+\log S_b\!\left(\frac{y+x}{\pi b}\right)
+\log S_b\!\left(\frac{y-x}{\pi b}\right)\\
&-2\log S_b\!\left(\frac{2y-\pi+x}{2\pi b}+\frac b2\right)
 -2\log S_b\!\left(\frac{2y-\pi-x}{2\pi b}+\frac b2\right)\Bigg].
\end{aligned}
\]
Here we use $\log S_b$ defined by the integral
in Section~2.3.
Define
\[
\begin{aligned}
V_r^+(x,\mathbf y)&=-\pi\left(1-\frac2r\right)x+\sum\limits_{i=1}^n\tau_b(x,y_i),\\
V_r^-(x,\mathbf y)&=\pi\left(1-\frac2r\right)x+\sum\limits_{i=1}^n\tau_b(x,y_i).
\end{aligned}
\]
Then
\[
\begin{aligned}
\mathrm{TV}_r(M_n)
={}&\frac{(-1)^{\frac{(r+1)(n+1)}2}\,2^{n-1}}{r^n}
\sin^{n-1}\frac{2\pi}{r}
\sum\limits_{\substack{a\in2\mathbb Z,\ \mathbf{k}\in\mathbb Z^n\\(x,\mathbf{y})\in D_S}}
\sum\limits_{\sigma\in\{+,-\}}e^{\frac{r}{2\pi\mathbf i}V_r^\sigma(x,\mathbf{y})}
+O\!\left(e^{\frac{r}{2\pi}(\operatorname{Vol}(M_n)-\varepsilon)}\right).
\end{aligned}
\]
Here $\varepsilon>0$ is independent of $n,r$, and the implied constant depends only on $n$.
\end{keylemma}

\begin{proof}
By Lemma~\ref{lem:hyperideal-colors}, it suffices to consider even colors
$a$ with $r-2<3a\le 2(r-2)$. Apply the integer sum \eqref{eq:reflected-integer}
of Theorem~\ref{thm:hyperideal-finite-sums} with $b^2=\frac{2}{r}$,
$i=1$ and $\mathbf a=(a,a,a,a,a,a)$. The parameters omitting $a_1$ are
$T_3=T_4=\frac{3a}{2}$ and $Q_3=2a$, so the theorem gives
\[
 \left|\begin{matrix}a&a&a\\a&a&a\end{matrix}\right|
 =\sum\limits_{k=\frac{3a}{2}}^{2a}U_b(F_1(\mathbf a);k),
 \qquad F_1(\mathbf a)=(r-2-a,a,a,a,a,a).
\]
The terms with $k>r-2$ vanish. Writing the original index as $a+k$ gives
$\frac a2\le k\le\min\{a,r-2-a\}$. Set $y=\frac{2\pi}{r}(k+1)$.
Evaluating $U_b(F_1(\mathbf a);a+k)$ by
\eqref{eq:root-kernel} of Proposition~\ref{prop:double-sine-sixj} gives
\begin{equation}\label{eq:one-edge-sb-summand}
\left|\begin{matrix}a&a&a\\a&a&a\end{matrix}\right|
=
\frac{2(-1)^{\frac{r+1}{2}}}{r}\sin\frac{2\pi}{r}
\sum\limits_{k=\frac a2}^{\min\{a,r-2-a\}}
\frac{S_b\!\left(\frac{y}{\pi b}\right)^2
      S_b\!\left(\frac{y+x}{\pi b}\right)S_b\!\left(\frac{y-x}{\pi b}\right)}
     {S_b\!\left(\frac{2y-\pi+x}{2\pi b}+\frac b2\right)^2
      S_b\!\left(\frac{2y-\pi-x}{2\pi b}+\frac b2\right)^2}.
\end{equation}
The product inside the sum in \eqref{eq:one-edge-sb-summand} equals
$\exp\!\left(\frac{r}{2\pi\mathbf i}\tau_b(x,y)\right)$.
At real summation points outside $D_S$, use $\log|S_b|$ in the two
terms with coefficient $-2$ in $\tau_b$.

Apply Proposition~\ref{prop:double-sine-asymptotics}(1) to each logarithm,
shifting $t$ by $-\frac{\pi b^2}{2}$ for the terms with arguments $y,y+x,y-x$.
Each scaled logarithm tends to $\frac12 L(t)$. Thus, as $b\to0$,
$\tau_b(x,y)$ converges locally uniformly for $\frac{\pi+|x|}{2}<y<\pi-|x|$ to
\[
\tau(x,y)=L(y)+\frac{L(y+x)+L(y-x)}{2}
-L\!\left(y-\frac{\pi+x}{2}\right)
-L\!\left(y-\frac{\pi-x}{2}\right).
\]
On $D_S$, the identity $L(t)=-2\mathbf i\Lambda(t)$ gives
\begin{equation}\label{eq:one-edge-classical-im}
\operatorname{Im}\tau(x,y)
=
2\Lambda\!\left(y-\frac{\pi+x}{2}\right)
+2\Lambda\!\left(y-\frac{\pi-x}{2}\right)
-2\Lambda(y)-\Lambda(y+x)-\Lambda(y-x).
\end{equation}
We use the right-hand side of \eqref{eq:one-edge-classical-im}
for the real exponent throughout $D_{H,r}$.

We claim that, uniformly over the terms of the sum in $D_{H,r}$,
\begin{equation}\label{eq:one-edge-log-estimate}
\operatorname{Im}\tau_b(x,y)=\operatorname{Im}\tau(x,y)
+O\!\left(\frac{\log r}{r}\right).
\end{equation}
By Proposition~\ref{prop:double-sine-asymptotics}(2), with $b^2=\frac2r$,
\begin{equation}\label{eq:centered-double-sine-bound}
\left|
\frac{2\pi\mathbf i}{r}\log S_b\!\left(\frac{t}{\pi b}+\frac b2\right)
-\frac12 L(t)
\right|\le\frac{B}{r},
\qquad 0\le t\le\pi.
\end{equation}
For $s=y,y+x,y-x$, apply \eqref{eq:centered-double-sine-bound} at
$t=s-\frac{\pi}{r}\in[0,\pi]$. Replacing $\Lambda(t)$ by $\Lambda(s)$
costs $O(\frac{\log r}{r})$ by Lemma~\ref{lem:lobachevsky-continuity}.
For the two remaining arguments $t=y-\frac{\pi\pm x}{2}$, apply
\eqref{eq:centered-double-sine-bound} directly whenever $t\ge0$.

Only $t=y-\frac{\pi-x}{2}$ can be negative. In that case
$-\frac{\pi}{3}<t<0$, so $t+\pi\in(0,\pi)$. Put $z=\frac{t}{\pi b}+\frac b2$.
The shift $t\mapsto t+\pi$ sends $z$ to $z+b^{-1}$, and the functional equation gives
\[
S_b(z+b^{-1})=2\cos\frac{rt}{2}\,S_b(z).
\]
Since $\frac{rt}{2\pi}=k+\frac a2-\frac{r-3}{2}\in\mathbb Z$ at these summation points,
we have $\log|S_b(z)|=\log S_b(z+b^{-1})-\log 2$.
Taking imaginary parts in \eqref{eq:centered-double-sine-bound} at $t+\pi$
now gives $\frac{2\pi}{r}\log|S_b(z)|=-\Lambda(t)+O(\frac1r)$,
since $\Lambda(t+\pi)=\Lambda(t)$. Combining the five estimates proves
\eqref{eq:one-edge-log-estimate}.

Let $v_8=8\Lambda(\frac{\pi}{4})$ be the volume of the regular ideal
octahedron. For $(x,\mathbf{y})\in D_S$ and $y=y_i$, we have
$\tau(-x,y)=\tau(x,y)$. For $x>0$,
$$
\operatorname{Im}\tau'_x(x,y)
=
\log\left(
\frac{\sin(y-\frac{\pi+x}{2})}{\sin(y-\frac{\pi-x}{2})}
\frac{\sin(y+x)}{\sin(y-x)}
\right)\le0.
$$
It follows that
$$
\operatorname{Im}\tau(x,y)\le\operatorname{Im}\tau(0,y).
$$
Since $y-\frac{\pi}{2}$ and $\pi-y$ have average $\frac{\pi}{4}$,
the concavity of $\Lambda$ on $[0,\frac{\pi}{2}]$ gives
\begin{equation}\label{eq:one-edge-octahedral-bound}
\operatorname{Im}\tau(0,y)
=4\left[\Lambda\!\left(y-\frac{\pi}{2}\right)+\Lambda(\pi-y)\right]
\le8\Lambda\!\left(\frac{\pi}{4}\right)=v_8.
\end{equation}

For $(x,\mathbf{y})\in D_{H,r}$, membership in $D_{H,r}\setminus D_S$
is equivalent to the existence of $i$ such that
$$
\left(x<0\ \text{ and }\ y_i\le\frac{\pi-x}{2}\right)
\qquad\text{or}\qquad y_i=\pi-|x|.
$$
In the first case,
$$
-\frac{\pi}{3}<y_i-\frac{\pi-x}{2}\le0,
\qquad
0<y_i+x\le\frac{\pi}{2}.
$$
Thus $2\Lambda(y_i-\frac{\pi-x}{2})$ and $-\Lambda(y_i+x)$ are
nonpositive. Bounding the remaining terms in
\eqref{eq:one-edge-classical-im} gives
$$
\operatorname{Im}\tau(x,y_i)\le5\Lambda\!\left(\frac{\pi}{6}\right)<2.77.
$$

In the second case, by symmetry we may take $x\ge0$ and $y_i=\pi-x$.
This includes $3a=2(r-2)$, where $k_i=\frac{a}{2}$ and $y_i=\pi-x$ for every $i$.
We have
$$
\begin{aligned}
\operatorname{Im}\tau(x,\pi-x)
&=2\Lambda\!\left(\frac{\pi-x}{2}\right)
 +2\Lambda\!\left(\frac{\pi-3x}{2}\right)
 +2\Lambda(x)+\Lambda(2x)\\
&\le4\Lambda\!\left(\frac{\pi}{6}\right)+x\log2+\Lambda(2x),
\end{aligned}
$$
where we used
$\Lambda(\frac{\pi}{2}-t)=\int_0^t\log(2\cos u)\,du\le t\log2$.
Since $\bigl(x\log2+\Lambda(2x)\bigr)'=\log2-2\log(2\sin 2x)$ for
$0<x\le\frac{\pi}{3}$, the maximum on $[0,\frac{\pi}{3}]$ occurs at $x=\frac{\pi}{8}$. Hence
$$
\operatorname{Im}\tau(x,\pi-x)
\le4\Lambda\!\left(\frac{\pi}{6}\right)
+\frac{\pi}{8}\log2+\Lambda\!\left(\frac{\pi}{4}\right)<2.77.
$$
Together with \eqref{eq:one-edge-octahedral-bound}, this proves that
for each omitted term,
\begin{equation}\label{eq:one-edge-exterior-bound}
\sum\limits_{i=1}^{n}\operatorname{Im}\tau(x,y_i)
\le(n-1)v_8+2.77.
\end{equation}

To compare \eqref{eq:one-edge-exterior-bound} with the volume,
substitute $u=\frac n2t$ in Lemma~\ref{lem:one-edge-volume}(1).
Since $n\ge2$ and $\frac{\cos t}{2\cos t-1}$ is increasing on
$[0,\frac{\pi}{6}]$, we have
$$
\begin{aligned}
nv_8-\operatorname{Vol}(M_n)
&=6\int_0^{\frac{\pi}{6}}
\operatorname{arcosh}\!\left(
\frac{\cos\frac{2u}{n}}{2\cos\frac{2u}{n}-1}
\right)\,du\\
&\le6\int_0^{\frac{\pi}{6}}
\operatorname{arcosh}\!\left(\frac{\cos u}{2\cos u-1}\right)\,du
=2v_8-\operatorname{Vol}(M_2).
\end{aligned}
$$
A numerical calculation gives $\operatorname{Vol}(M_2)-v_8\approx 2.788128$.
Therefore
$$
\operatorname{Vol}(M_n)-\bigl((n-1)v_8+2.77\bigr)
\ge\operatorname{Vol}(M_2)-v_8-2.77>0.01.
$$
For fixed $n$, the edge weight and the number of summands grow at most
polynomially in $r$. It follows from \eqref{eq:one-edge-log-estimate},
\eqref{eq:one-edge-exterior-bound}, and Lemma~\ref{lem:hyperideal-colors}
that
\begin{equation}\label{eq:one-edge-even-core}
\begin{aligned}
\mathrm{TV}_r(M_n)
={}&\left(\frac{2(-1)^{\frac{r+1}{2}}}{r}\sin\frac{2\pi}{r}\right)^n
\sum\limits_{\substack{a\in2\mathbb Z,\ \mathbf{k}\in\mathbb Z^n\\
                 (x,\mathbf{y})\in D_S}}
[a+1]\exp\!\left(\frac{r}{2\pi\mathbf i}\sum\limits_{i=1}^n \tau_b(x,y_i)\right)\\
&{}+O\!\left(e^{\frac r{2\pi}(\operatorname{Vol}(M_n)-\varepsilon)}\right)
\end{aligned}
\end{equation}
for some $\varepsilon>0$ independent of $n,r$. Indeed, for even $a$ and
integer $k_i$, the inequalities defining $D_S$ give
$$
r-2<3a<2(r-2),\qquad
2k_i\ge\max\{a,r-2-a\},\qquad
k_i<\min\{a,r-2-a\},
$$
so all indices in \eqref{eq:one-edge-even-core} belong to the original summation range.
For integer $a$, the definition of $x$ gives
\begin{equation}\label{eq:one-edge-signed-weight}
(-1)^a[a+1]
=
\frac{(-1)^{\frac{r+1}{2}}}{2\sin\frac{2\pi}{r}}
\sum\limits_{\sigma\in\{+,-\}}
\exp\!\left(\sigma\mathbf i\left(\frac r2-1\right)x\right).
\end{equation}
Since $a$ is even in \eqref{eq:one-edge-even-core}, substituting
\eqref{eq:one-edge-signed-weight} for $[a+1]$ gives the stated sum with potential functions $V_r^+$ and $V_r^-$.
\end{proof}

Since $\tau_b\to\tau$ as $b\to0$, the classical potential functions on $D_S$ are
\[
\begin{aligned}
V^+(x,\mathbf y)&=\lim\limits_{r\to\infty}V_r^+(x,\mathbf y)=-\pi x+\sum\limits_{i=1}^n\tau(x,y_i),\\
V^-(x,\mathbf y)&=\lim\limits_{r\to\infty}V_r^-(x,\mathbf y)=\pi x+\sum\limits_{i=1}^n\tau(x,y_i).
\end{aligned}
\]

\subsection{Reflection and summation over all integers}\label{subsec:reflection-extension}

When Poisson summation is applied directly to the sum over even $a$, the
potential function of the central Fourier term gives the growth rate of its
exponential factor on the scale $\frac{r}{2\pi}$. At $x=0$ and
$y_1=\cdots=y_n=\frac{3\pi}{4}$, \eqref{eq:one-edge-octahedral-bound} and
Lemma~\ref{lem:one-edge-volume}(1) give
\[
\sum\limits_{i=1}^n\operatorname{Im}\tau(x,y_i)
=n\operatorname{Im}\tau\!\left(0,\frac{3\pi}{4}\right)
=8n\Lambda\!\left(\frac{\pi}{4}\right)=nv_8>\operatorname{Vol}(M_n).
\]
Here both arguments of $\Lambda$ in \eqref{eq:one-edge-octahedral-bound}
equal $\frac{\pi}{4}$.
This excess over $\operatorname{Vol}(M_n)$ is the main obstacle in the traditional approach.

In the specialization of \eqref{eq:reflected-integer} to equal colors,
the double sine factors in \eqref{eq:one-edge-sb-summand} are permuted by
$x\mapsto-x$. Hence the potentials of Key Lemma~\ref{lem:symmetric-region} satisfy
\begin{equation}\label{eq:quantum-potential-symmetry}
V_r^+(-x,\mathbf y)=V_r^-(x,\mathbf y).
\end{equation}
The next lemma uses this symmetry to extend the sum from even $a$ to all
integers $a$, thereby avoiding the obstacle above.

\begin{keylemma}\label{lem:reflection-extension}
Let $r\ge3$ be odd. Then
\[
\sum\limits_{\substack{a\in2\mathbb Z,\ \mathbf k\in\mathbb Z^n\\(x,\mathbf y)\in D_S}}
\sum\limits_{\sigma\in\{+,-\}}e^{\frac{r}{2\pi\mathbf i}V_r^\sigma(x,\mathbf y)}
=\frac12
\sum\limits_{\substack{(a,\mathbf k)\in\mathbb Z^{n+1}\\(x,\mathbf y)\in D_S}}
\sum\limits_{\sigma\in\{+,-\}}e^{\frac{r}{2\pi\mathbf i}V_r^\sigma(x,\mathbf y)}.
\]
\end{keylemma}

\begin{proof}
The involution $a\mapsto r-2-a$ sends $x$ to $-x$ and preserves $D_S$.
By \eqref{eq:quantum-potential-symmetry} it interchanges the two exponential
terms. Since $r$ is odd, it also interchanges even and odd $a$, so the full
sum is twice the even sum.
\end{proof}

Combining Key Lemmas~\ref{lem:symmetric-region} and~\ref{lem:reflection-extension} gives
\begin{equation}\label{eq:one-edge-integer-sum}
\mathrm{TV}_r(M_n)
=\frac{(-1)^{\frac{(r+1)(n+1)}2}\,2^{n-2}}{r^n}
\sin^{n-1}\frac{2\pi}{r}
\sum\limits_{\substack{(a,\mathbf k)\in\mathbb Z^{n+1}\\(x,\mathbf y)\in D_S}}
\sum\limits_{\sigma\in\{+,-\}} e^{\frac{r}{2\pi\mathbf i}V_r^\sigma(x,\mathbf y)}
+O\!\left(e^{\frac{r}{2\pi}(\operatorname{Vol}(M_n)-\varepsilon)}\right).
\end{equation}

\subsection{Poisson summation}

We first recall the Poisson summation formula
\cite[Theorem~3.1]{SteinShakarchi}. For $f\in C_c^\infty(\mathbb R^{n+1})$,
let
$$
\widehat f(\mathbf{b})
=\int_{\mathbb R^{n+1}}f(\mathbf{z})\,
e^{-2\pi\mathbf i\,\mathbf{b}\cdot\mathbf{z}}\,d\mathbf{z}.
$$
Then
\begin{equation}\label{eq:poisson-summation}
\sum\limits_{\mathbf{a}\in\mathbb Z^{n+1}}f(\mathbf{a})
=\sum\limits_{\mathbf{b}\in\mathbb Z^{n+1}}\widehat f(\mathbf{b}).
\end{equation}

To apply Poisson summation \eqref{eq:poisson-summation} to the sum in
\eqref{eq:one-edge-integer-sum}, let
$$
D_{S,\delta}
=\left\{(x,\mathbf{y})\in D_S:
\operatorname{dist}\bigl((x,\mathbf{y}),\partial D_S\bigr)>\delta
\right\}.
$$
For sufficiently small $\delta>0$, choose a smooth bump function
$\psi(x,\mathbf{y})$ with values in $[0,1]$ such that
$\psi(-x,\mathbf{y})=\psi(x,\mathbf{y})$ and
$$
\psi(x,\mathbf{y})=
\begin{cases}
1, & (x,\mathbf{y})\in D_{S,2\delta},\\
0, & (x,\mathbf{y})\notin D_{S,\delta}.
\end{cases}
$$

\begin{lemma}\label{lem:one-edge-poisson}
As $r\to\infty$ through odd integers,
\begin{equation}\label{eq:one-edge-poisson}
\begin{aligned}
&\mathrm{TV}_r(M_n)
={}
\frac{(-1)^{\frac{(r+1)(n+1)}{2}}\,r}{8\pi^{n+1}}
\sin^{n-1}\frac{2\pi}{r}
\sum\limits_{\sigma\in\{+,-\}}
\sum\limits_{(p,\mathbf{q})\in\mathbb Z^{n+1}}(-1)^p\\
&\times
\int_{D_S}\psi(x,\mathbf{y})\,
\exp\!\left[
\frac{r}{2\pi\mathbf i}
\left(V_r^\sigma(x,\mathbf{y})
+2\pi p x+2\pi\mathbf{q}\cdot\mathbf{y}\right)
\right]\,dx\,d\mathbf{y}+
O\!\left(e^{\frac{r}{2\pi}(\operatorname{Vol}(M_n)-\varepsilon)}\right)
\end{aligned}
\end{equation}
for some $\varepsilon>0$ independent of $n$ and $r$.
The implied constant depends only on $n$.
\end{lemma}

\begin{proof}
Similarly, by \eqref{eq:one-edge-classical-im} and the boundary estimates
in the proof of Key Lemma~\ref{lem:symmetric-region}, for sufficiently small
$\delta>0$ we have
$$
\sup\limits_{(x,\mathbf{y})\in D_S\setminus D_{S,2\delta}}
\sum\limits_{i=1}^{n}\operatorname{Im}\tau(x,y_i)
\le\operatorname{Vol}(M_n)-\varepsilon'
$$
for some $\varepsilon'>0$ independent of $n$ and $r$.
By \eqref{eq:centered-double-sine-bound} and
Lemma~\ref{lem:lobachevsky-continuity}, uniformly on $D_S$,
$$
\operatorname{Im}V_r^\pm(x,\mathbf{y})
=\sum\limits_{i=1}^{n}\operatorname{Im}\tau(x,y_i)
+O\!\left(\frac{n\log r}{r}\right).
$$
Since $\psi=1$ on $D_{S,2\delta}$ and $O(r^{n+1})$ integer tuples
$(a,\mathbf k)$ satisfy $(x,\mathbf y)\in D_S$, inserting $\psi(x,\mathbf{y})$
in \eqref{eq:one-edge-integer-sum} introduces an error of order
$O\!\left(e^{\frac{r}{2\pi}(\operatorname{Vol}(M_n)-\varepsilon)}\right)$
for some $\varepsilon>0$.

For each sign, apply \eqref{eq:poisson-summation}
to the sum in $a,\mathbf{k}$. In the $(p,\mathbf{q})$-th
Fourier integral, make the change of variables
$$
a=\frac{r}{2\pi}(x+\pi)-1,\qquad
k_i=\frac{r}{2\pi}y_i-1,\qquad 1\le i\le n.
$$
Then
$$
da\,d\mathbf{k}
=\left(\frac{r}{2\pi}\right)^{n+1}dx\,d\mathbf{y},
\qquad
e^{-2\pi\mathbf i(pa+\mathbf{q}\cdot\mathbf{k})}
=(-1)^p e^{-\mathbf i r(px+\mathbf{q}\cdot\mathbf{y})},
$$
since $r$ is odd. For either sign, the $(p,\mathbf{q})$-th Fourier coefficient
equals
$$
(-1)^p\left(\frac{r}{2\pi}\right)^{n+1}
\int_{D_S}\psi(x,\mathbf{y})\,
\exp\!\left[
\frac{r}{2\pi\mathbf i}
\left(V_r^\pm(x,\mathbf{y})
+2\pi p x+2\pi\mathbf{q}\cdot\mathbf{y}\right)
\right]\,dx\,d\mathbf{y}.
$$
Multiplying these coefficients by the prefactor in \eqref{eq:one-edge-integer-sum}
and summing over both signs and all $(p,\mathbf q)$ gives \eqref{eq:one-edge-poisson}.
\end{proof}

In Section~5 we estimate the Fourier terms in \eqref{eq:one-edge-poisson}.
For each $(p,\mathbf q)$, write
$$
f_r^\pm(p,\mathbf{q})
=(-1)^p\int_{D_S}\psi(x,\mathbf{y})\,
\exp\!\left[
\frac{r}{2\pi\mathbf i}
\left(V_r^\pm(x,\mathbf{y})+2\pi p x
+2\pi\mathbf{q}\cdot\mathbf{y}\right)
\right]\,dx\,d\mathbf{y}.
$$

Since $D_S$ and $\psi$ are invariant under $x\mapsto-x$,
the substitution $x\mapsto-x$ and \eqref{eq:quantum-potential-symmetry} give
$$
f_r^+(p,\mathbf{q})=f_r^-(-p,\mathbf{q}).
$$

The same change of variables in $f_r^\pm(\pm1-p,\mathbf{q})$,
together with the identity
$$
V_r^\pm(-x,\mathbf{y})-2\pi(\pm1-p)x
=V_r^\pm(x,\mathbf{y})+2\pi p x\mp\frac{4\pi x}{r},
$$
gives
\begin{equation}\label{eq:fourier-index-reflection}
f_r^\pm(\pm1-p,\mathbf{q})
=-(-1)^p\int_{D_S}e^{\pm2\mathbf i x}\,\psi(x,\mathbf{y})\,
\exp\!\left[
\frac{r}{2\pi\mathbf i}
\left(V_r^\pm(x,\mathbf{y})+2\pi p x
+2\pi\mathbf{q}\cdot\mathbf{y}\right)
\right]\,dx\,d\mathbf{y}.
\end{equation}

\section{\texorpdfstring{Asymptotic analysis}{Asymptotic analysis}}

\subsection{Geometric critical points and Reidemeister torsion}

Recall that the two classical potentials are
$$
V^+(x,\mathbf{y})=-\pi x+\sum\limits_{i=1}^{n}\tau(x,y_i),
\qquad
V^-(x,\mathbf{y})=\pi x+\sum\limits_{i=1}^{n}\tau(x,y_i).
$$
Since $\tau(-x,y)=\tau(x,y)$, we have
\begin{equation}\label{eq:classical-potential-symmetry}
V^-(x,\mathbf{y})=V^+(-x,\mathbf{y}).
\end{equation}
Thus, if $(x_0,\mathbf{y}_0)$ is a critical point of $V^+$,
then $(-x_0,\mathbf{y}_0)$ is a critical point of $V^-$ with the same
critical value. Differentiating \eqref{eq:classical-potential-symmetry}
twice also gives
$$
\det\operatorname{Hess}V^-(-x_0,\mathbf{y}_0)
=\det\operatorname{Hess}V^+(x_0,\mathbf{y}_0).
$$

\begin{proposition}\label{prop:geometric-critical-points}
Let $\ell_n$ be the length of the unique edge of $M_n$, so that
$$
\cosh\ell_n
=\frac{\cos\frac{\pi}{3n}}
{2\cos\frac{\pi}{3n}-1}.
$$
Let $y^{\mathrm{root}}(n)\in(\frac{\pi}{2},\pi)$ be the unique solution of
$$
4\cosh\ell_n\,\sin^2 y^{\mathrm{root}}(n)=1+\cosh\ell_n.
$$
Then
$$
\mathbf{c}^+=(\mathbf i\ell_n,y^{\mathrm{root}}(n),\ldots,y^{\mathrm{root}}(n)),
\qquad
\mathbf{c}^-=(-\mathbf i\ell_n,y^{\mathrm{root}}(n),\ldots,y^{\mathrm{root}}(n))
$$
are critical points of $V^+$ and $V^-$, respectively, and
$$
V^+(\mathbf{c}^+)=V^-(\mathbf{c}^-)=\mathbf i\,\operatorname{Vol}(M_n).
$$
\end{proposition}

\begin{proof}
By \eqref{eq:classical-potential-symmetry}, it suffices to consider
$V^+$. For $0<\theta\le\frac{\pi}{3n}$, let $\ell(\theta)$ be the
common length of the internal edges of the regular truncated tetrahedron
with dihedral angle $\theta$. By \cite[Eq.~(3.1)]{AtkinsonRafalski},
$$
\cosh\ell(\theta)=\frac{\cos\theta}{2\cos\theta-1}.
$$
Set $\ell(0)=0$. For $0\le\theta\le\frac{\pi}{3n}$, let
$y(\theta)\in(\frac{\pi}{2},\pi)$ be determined by
$$
4\cosh\ell(\theta)\,\sin^2 y(\theta)=1+\cosh\ell(\theta).
$$
Write $\ell=\ell(\theta)$ and $y=y(\theta)$. Using
$L'(z)=2\mathbf i\log(2\sin z)$ in the definition of $\tau$, we obtain
$$
\tau'_y(\mathbf i\ell,y)
=2\mathbf i\log\!\left(
\frac{\sin y\,\sqrt{\sin^2 y+\sinh^2\ell}}
{\cos^2 y+\sinh^2\frac{\ell}{2}}
\right).
$$
The relation $4\cosh\ell\,\sin^2 y=1+\cosh\ell$ gives
$$
\sin y\,\sqrt{\sin^2 y+\sinh^2\ell}
=\cos^2 y+\sinh^2\frac{\ell}{2},
$$
so $\tau'_y(\mathbf i\ell,y)=0$. Similarly,
$$
\tau'_x(\mathbf i\ell,y)
=2\arctan\!\left(-\cot y\,\tanh\ell\right)
+2\arctan\!\left(-\tan y\,\tanh\frac{\ell}{2}\right).
$$
The defining equations for $\ell$ and $y$ give
$$
-\cot y\,\tanh\ell=\tan\theta,
\qquad
-\tan y\,\tanh\frac{\ell}{2}=\tan\frac{\theta}{2}.
$$
Hence $\tau'_x(\mathbf i\ell,y)=3\theta$. At $\theta=\frac{\pi}{3n}$,
we have $\ell=\ell_n$ and $y=y^{\mathrm{root}}(n)$, so
$$
(V^+)'_{y_i}(\mathbf{c}^+)=0
\quad (1\le i\le n),
\qquad
(V^+)'_x(\mathbf{c}^+)
=-\pi+3n\frac{\pi}{3n}=0.
$$
Thus $\mathbf{c}^+$ is a critical point of $V^+$.

To compute the critical value, apply the chain rule with
$\tau'_y(\mathbf i\ell,y)=0$ and $\tau'_x(\mathbf i\ell,y)=3\theta$. This gives
$$
\left[
\tau\bigl(\mathbf i\ell(\theta),y(\theta)\bigr)
-3\mathbf i\theta\,\ell(\theta)
\right]'_{\theta}
=-3\mathbf i\,\ell(\theta).
$$
At $\theta=0$, we have $y(0)=\frac{3\pi}{4}$ and
$$
\tau\!\left(0,\frac{3\pi}{4}\right)
=2L\!\left(\frac{3\pi}{4}\right)-2L\!\left(\frac{\pi}{4}\right)
=8\mathbf i\Lambda\!\left(\frac{\pi}{4}\right)
=\mathbf i v_8.
$$
Integrating from $0$ to $\frac{\pi}{3n}$ and applying
Lemma~\ref{lem:one-edge-volume}(1), we conclude that
$$
\begin{aligned}
V^+(\mathbf{c}^+)
&=-\pi\mathbf i\ell_n+n\tau(\mathbf i\ell_n,y^{\mathrm{root}}(n))\\
&=\mathbf i n\left[
v_8-3\int_0^{\frac{\pi}{3n}}
\operatorname{arcosh}\!\left(\frac{\cos\theta}{2\cos\theta-1}\right)
\,d\theta
\right]
=\mathbf i\,\operatorname{Vol}(M_n).
\end{aligned}
$$
This proves the proposition.
\end{proof}
\begin{proposition}\label{prop:hessian-torsion}
Let $DM_n$ be the double of $M_n$ along its totally geodesic boundary.
Write $\rho\colon\pi_1(DM_n)\to\mathrm{PSL}(2,\mathbb C)$ for its
hyperbolic holonomy representation. At the critical points $\mathbf{c}^+$ and
$\mathbf{c}^-$ of Proposition~\ref{prop:geometric-critical-points}, the Hessians
with respect to $x,y_1,\ldots,y_n$ satisfy
$$
\begin{aligned}
&\det\!\left(\mathbf i\,\operatorname{Hess}V^+(\mathbf{c}^+)\right)
=\det\!\left(\mathbf i\,\operatorname{Hess}V^-(\mathbf{c}^-)\right)\\
=&4\sinh^2\ell_n
\left(
\frac{\cosh\ell_n}
{(1+\cosh\ell_n)(2\cosh\ell_n-1)}
\right)^{2n}
\left|\operatorname{Tor}(DM_n;\operatorname{Ad}\rho)\right|,
\end{aligned}
$$
where $\operatorname{Tor}(DM_n;\operatorname{Ad}\rho)$ denotes the
adjoint Reidemeister torsion of $DM_n$. In particular, both critical
points are nondegenerate.
\end{proposition}

\begin{proof}
By \eqref{eq:classical-potential-symmetry}, it suffices to compute the
Hessian of $V^+$ at $\mathbf{c}^+$. In what follows, the second derivatives of
$\tau$ are evaluated at $(\mathbf i\ell_n,y^{\mathrm{root}}(n))$. Since
$$
V^+(x,\mathbf{y})=-\pi x+\sum\limits_{i=1}^n \tau(x,y_i),
$$
we have
$$
\operatorname{Hess}V^+(\mathbf{c}^+)
=
\begin{pmatrix}
n\tau''_{xx} & \tau''_{xy} & \cdots & \tau''_{xy}\\
\tau''_{xy} & \tau''_{yy} & & 0\\
\vdots & & \ddots & \\
\tau''_{xy} & 0 & & \tau''_{yy}
\end{pmatrix}.
$$

Write $c=\cosh\ell_n$. By the expression for $\tau'_y$ in the proof of
Proposition~\ref{prop:geometric-critical-points}, for $y$ near $y^{\mathrm{root}}(n)$
we have
$$
\tau'_y(\mathbf i\ell_n,y)
=\mathbf i\log\!\left(
\frac{\sin^2 y\,(\sin^2 y+c^2-1)}
{\left(\frac{1+c}{2}-\sin^2 y\right)^2}
\right).
$$
Differentiating with respect to $y$ and using
$\sin^2 y^{\mathrm{root}}(n)=\frac{1+c}{4c}$ and $\frac{\pi}{2}<y^{\mathrm{root}}(n)<\pi$, we obtain
\begin{equation}\label{eq:one-edge-hessian-yy}
\mathbf i \tau''_{yy}
=\frac{8c^2}{(2c-1)^2}\sqrt{\frac{3c-1}{c+1}}>0.
\end{equation}

For each $i=1,\ldots,n$, subtract $\frac{\tau''_{xy}}{\tau''_{yy}}$ times
row $i+1$ from the first row of $\operatorname{Hess}V^+(\mathbf{c}^+)$.
The first row becomes
$$
\left(n\left(\tau''_{xx}-\frac{(\tau''_{xy})^2}{\tau''_{yy}}\right),0,\ldots,0\right).
$$
Expanding along this row and multiplying by $\mathbf i^{n+1}$ gives
$$
\det\!\left(\mathbf i\,\operatorname{Hess}V^+(\mathbf{c}^+)\right)
=(\mathbf i \tau''_{yy})^n\,n\mathbf i
\left(\tau''_{xx}-\frac{(\tau''_{xy})^2}{\tau''_{yy}}\right).
$$

To evaluate the last factor, recall from the proof of
Proposition~\ref{prop:geometric-critical-points} that
$$
\tau'_y\bigl(\mathbf i\ell(\theta),y(\theta)\bigr)=0,
\qquad
\tau'_x\bigl(\mathbf i\ell(\theta),y(\theta)\bigr)=3\theta.
$$
Differentiating at $\theta_n=\frac{\pi}{3n}$, we obtain
$$
\begin{aligned}
\mathbf i\ell'(\theta_n)\tau''_{xy}+y'(\theta_n)\tau''_{yy}&=0,\\
\mathbf i\ell'(\theta_n)\tau''_{xx}+y'(\theta_n)\tau''_{xy}&=3.
\end{aligned}
$$
Eliminating $y'(\theta_n)$ gives
$$
\mathbf i\left(\tau''_{xx}-\frac{(\tau''_{xy})^2}{\tau''_{yy}}\right)
=\frac{3}{\ell'(\theta_n)}.
$$
Hence
\begin{equation}\label{eq:one-edge-hessian-determinant}
\det\!\left(\mathbf i\,\operatorname{Hess}V^+(\mathbf{c}^+)\right)
=\frac{3n}{\ell'(\theta_n)}(\mathbf i \tau''_{yy})^n.
\end{equation}

Each of the $n$ regular truncated tetrahedra has the same length Gram
matrix
$$
G=
\begin{pmatrix}
1&-c&-c&-c\\
-c&1&-c&-c\\
-c&-c&1&-c\\
-c&-c&-c&1
\end{pmatrix},
\qquad
-\det G=(3c-1)(1+c)^3.
$$
Thus \eqref{eq:one-edge-hessian-yy} can be written as
\begin{equation}\label{eq:one-edge-hessian-gram}
\mathbf i \tau''_{yy}
=8\left(\frac{c}{(1+c)(2c-1)}\right)^2\sqrt{-\det G}.
\end{equation}

The total cone angle around the unique edge is $6n\theta$, so its
derivative with respect to the edge length at $\ell_n$ is
$\frac{6n}{\ell'(\theta_n)}$. By Wong--Yang
\cite[Theorem~1.6(3)]{WongYangTorsion},
$$
\left|\operatorname{Tor}(DM_n;\operatorname{Ad}\rho)\right|
=
2^{3n-3}\frac{6n}{\ell'(\theta_n)}
\frac{(-\det G)^{\frac{n}{2}}}{\sinh^2\ell_n}.
$$
Substituting \eqref{eq:one-edge-hessian-gram} into
\eqref{eq:one-edge-hessian-determinant} and comparing with the torsion
formula, we obtain
$$
\det\!\left(\mathbf i\,\operatorname{Hess}V^+(\mathbf{c}^+)\right)
=
4\sinh^2\ell_n
\left(
\frac{\cosh\ell_n}
{(1+\cosh\ell_n)(2\cosh\ell_n-1)}
\right)^{2n}
\left|\operatorname{Tor}(DM_n;\operatorname{Ad}\rho)\right|.
$$
Since $\ell'(\theta_n)>0$, equations \eqref{eq:one-edge-hessian-yy}
and \eqref{eq:one-edge-hessian-determinant} show that the determinant
is positive. The same conclusion holds for $V^-$ at $\mathbf{c}^-$ by
\eqref{eq:classical-potential-symmetry}.
\end{proof}
\subsection{Asymptotics of the leading terms}

The four leading terms in \eqref{eq:one-edge-poisson} are
\begin{align*}
f_r^+(0,\mathbf{0})
&=\int_{D_S}\psi(x,\mathbf{y})\,
e^{\frac{r}{2\pi\mathbf i}V_r^+(x,\mathbf{y})}\,dx\,d\mathbf{y},\\
f_r^-(0,\mathbf{0})
&=\int_{D_S}\psi(x,\mathbf{y})\,
e^{\frac{r}{2\pi\mathbf i}V_r^-(x,\mathbf{y})}\,dx\,d\mathbf{y},\\
f_r^+(1,\mathbf{0})
&=-\int_{D_S}\psi(x,\mathbf{y})\,
e^{\frac{r}{2\pi\mathbf i}(V_r^+(x,\mathbf{y})+2\pi x)}
\,dx\,d\mathbf{y},\\
f_r^-(-1,\mathbf{0})
&=-\int_{D_S}\psi(x,\mathbf{y})\,
e^{\frac{r}{2\pi\mathbf i}(V_r^-(x,\mathbf{y})-2\pi x)}
\,dx\,d\mathbf{y}.
\end{align*}

The identity $f_r^+(p,\mathbf q)=f_r^-(-p,\mathbf q)$ gives
$$
f_r^-(0,\mathbf{0})=f_r^+(0,\mathbf{0}),\qquad
f_r^-(-1,\mathbf{0})=f_r^+(1,\mathbf{0}).
$$
Taking the upper signs, $p=0$ and $\mathbf{q}=\mathbf{0}$ in
\eqref{eq:fourier-index-reflection} gives
$$
f_r^+(1,\mathbf{0})
=-\int_{D_S}e^{2\mathbf i x}\psi(x,\mathbf{y})\,
e^{\frac{r}{2\pi\mathbf i}V_r^+(x,\mathbf{y})}\,dx\,d\mathbf{y}.
$$
Hence
\begin{equation}\label{eq:one-edge-four-term-sum}
\begin{aligned}
&f_r^+(0,\mathbf{0})+f_r^+(1,\mathbf{0})
 +f_r^-(0,\mathbf{0})+f_r^-(-1,\mathbf{0})\\
=&2\int_{D_S}(1-e^{2\mathbf i x})\psi(x,\mathbf{y})\,
e^{\frac{r}{2\pi\mathbf i}V_r^+(x,\mathbf{y})}\,dx\,d\mathbf{y}.
\end{aligned}
\end{equation}

To estimate \eqref{eq:one-edge-four-term-sum}, we first record the
following properties of $\operatorname{Im}\tau(x+\mathbf i t,y)$.

\begin{lemma}\label{lem:one-edge-shifted-potential}
Let $x,y\in\mathbb R$ satisfy $\frac{\pi+|x|}{2}<y<\pi-|x|$, and let $0\le t\le\ell_n$.
\begin{enumerate}
\renewcommand{\labelenumi}{(\arabic{enumi})}
\item For fixed $t$ and $y$, $\operatorname{Im}\tau(x+\mathbf i t,y)$
is even in $x$ and strictly decreasing for $x>0$.
\item For fixed $t$, the function $y\mapsto\operatorname{Im}\tau(\mathbf i t,y)$
on $(\frac{\pi}{2},\pi)$ attains its unique maximum at the solution of
$$
4\cosh t\,\sin^2y=1+\cosh t.
$$
In particular, when $t=\ell_n$, the maximum is attained at $y=y^{\mathrm{root}}(n)$.
\item For fixed $x,y$ with $|x|\le\frac{\pi}{4}$, the function
$t\mapsto\operatorname{Im}\tau(x+\mathbf i t,y)$ is convex on $[0,\ell_n]$.
\end{enumerate}
\end{lemma}

\begin{proof}
Differentiating the expression for $\tau$ in the proof of
Key Lemma~\ref{lem:symmetric-region} gives
$$
2\operatorname{Im}\tau'_x(x+\mathbf i t,y)
=\log\frac{\cosh t+\cos(2y-x)}{\cosh t+\cos(2y+x)}
+\log\frac{\cosh(2t)-\cos(2y+2x)}{\cosh(2t)-\cos(2y-2x)}.
$$
The right-hand side is odd in $x$. For $x>0$, both ratios are less
than $1$, since $\sin(2y)<0$. This proves (1).

For (2), differentiating the expression for
$\operatorname{Im}\tau'_y(\mathbf i t,y)$ gives
$$
\operatorname{Im}\tau''_{yy}(\mathbf i t,y)
=\sin(2y)\left(
\frac{1}{\sin^2 y}
+\frac{1}{\sin^2 y+\sinh^2 t}
+\frac{2}{\frac{1+\cosh t}{2}-\sin^2 y}
\right)<0.
$$
By the proof of Proposition~\ref{prop:geometric-critical-points},
the first derivative vanishes when $4\cosh t\,\sin^2y=1+\cosh t$.
Strict concavity therefore gives (2).

For (3), write $z=x+\mathbf i t$. Differentiating twice gives
$$
\bigl[\operatorname{Im}\tau(x+\mathbf i t,y)\bigr]''_{tt}
=\operatorname{Re}\!\left[
\frac12\cot\!\left(y-\frac{\pi+z}{2}\right)
+\frac12\cot\!\left(y-\frac{\pi-z}{2}\right)
-\cot(y+z)-\cot(y-z)\right].
$$
The first two cotangents have positive real parts, since their
arguments have real parts in $(0,\frac{\pi}{2})$. Moreover,
$$
\operatorname{Re}\bigl(\cot(y+z)+\cot(y-z)\bigr)
=\frac{2\sin(2y)\bigl(\cos(2x)\cosh(2t)-\cos(2y)\bigr)}
{\bigl(\cosh(2t)-\cos(2y+2x)\bigr)
 \bigl(\cosh(2t)-\cos(2y-2x)\bigr)}.
$$
For $|x|\le\frac{\pi}{4}$, we have $\cos(2x)\ge0$ and
$\cos(2y)\le\cos(2x)$. Since $\sin(2y)<0$, the last expression
is nonpositive. This proves (3).
\end{proof}

We now apply the following saddle point approximation to the integral in
\eqref{eq:one-edge-four-term-sum}.

\begin{proposition}[{\cite[Proposition~5.2]{MWWY}}]
\label{prop:saddle-point-method}
Let $D$ be a region in $\mathbb C^n$ and let $f(\mathbf{z})$ and
$g(\mathbf{z})$ be holomorphic functions on $D$ independent of $r$.
Let $f_r(\mathbf{z})$ be a holomorphic function of the form
$$
f_r(\mathbf{z})
=f(\mathbf{z})+\frac{\upsilon_r(\mathbf{z})}{r^2}.
$$
Let $S$ be an embedded real $n$-dimensional closed disk in $D$ and
let $\mathbf{c}$ be a point on $S$. If
\begin{enumerate}
\renewcommand{\labelenumi}{(\roman{enumi})}
\item $\mathbf{c}$ is a critical point of $f$ in $D$,
\item $\operatorname{Re}(f)(\mathbf{c})>
      \operatorname{Re}(f)(\mathbf{z})$ for all
      $\mathbf{z}\in S\setminus\{\mathbf{c}\}$,
\item the domain
      $\{\mathbf{z}\in D\mid
      \operatorname{Re}f(\mathbf{z})<\operatorname{Re}f(\mathbf{c})\}$
      deformation retracts to $S\setminus\{(\mathbf{c})\}$,
\item the Hessian matrix $\operatorname{Hess}(f)(\mathbf{c})$
      of $f$ at $\mathbf{c}$ is non-singular,
\item $g(\mathbf{c})\ne0$, and
\item $|\upsilon_r(\mathbf{z})|$ is bounded from above by a
      constant independent of $r$ in $D$,
\end{enumerate}
then
$$
\int_S g(\mathbf{z})e^{r f_r(\mathbf{z})}\,d\mathbf{z}
=\left(\frac{2\pi}{r}\right)^{\frac n2}
\frac{g(\mathbf{c})}
{\sqrt{\det\!\left(-\operatorname{Hess}(f)(\mathbf{c})\right)}}
e^{r f(\mathbf{c})}\left(1+O\!\left(\frac1r\right)\right).
$$
\end{proposition}

\begin{proposition}\label{prop:one-edge-leading-terms}
For fixed $n\ge2$, as $r\to\infty$ through odd integers,
\begin{equation}\label{eq:one-edge-leading-asymptotic}
\begin{aligned}
&f_r^+(0,\mathbf{0})+f_r^+(1,\mathbf{0})
 +f_r^-(0,\mathbf{0})+f_r^-(-1,\mathbf{0})\\
&\qquad=
2\left(\frac{4\pi^2}{r}\right)^{\frac{n+1}{2}}
\frac{e^{\frac{r}{2\pi}\operatorname{Vol}(M_n)}}
{\sqrt{\left|\operatorname{Tor}(DM_n;\operatorname{Ad}\rho)\right|}}
\left(1+O\!\left(\frac1r\right)\right).
\end{aligned}
\end{equation}
\end{proposition}
\begin{proof}
To apply Proposition~\ref{prop:saddle-point-method} at $\mathbf{c}^+$, we first estimate
the contribution of $D_S\setminus D_{S,2\delta}$ to
\eqref{eq:one-edge-four-term-sum}. By the proof of
Lemma~\ref{lem:one-edge-poisson}, for sufficiently small $\delta>0$
and all sufficiently large $r$,
$$
\sup\limits_{(x,\mathbf{y})\in D_S\setminus D_{S,2\delta}}
\operatorname{Im}V_r^+(x,\mathbf{y})
\le\operatorname{Vol}(M_n)-\varepsilon'
$$
for some $\varepsilon'>0$. Since $|1-e^{2\mathbf i x}|\le2$ on the
real domain, $0\le\psi\le1$, and $\psi=1$ on $D_{S,2\delta}$,
\eqref{eq:one-edge-four-term-sum} equals
$$
2\int_{D_{S,2\delta}}(1-e^{2\mathbf i x})
e^{\frac{r}{2\pi\mathbf i}V_r^+(x,\mathbf{y})}
\,dx\,d\mathbf{y}
+O\!\left(e^{\frac{r}{2\pi}(\operatorname{Vol}(M_n)-\varepsilon')}\right).
$$

We translate $x$ to $x+\mathbf i\ell_n$, keeping $\mathbf{y}$
fixed. Denote the translated domain and the sides swept out by its
boundary by
$$
S_{\mathrm{top}}
=\left\{(x+\mathbf i\ell_n,\mathbf{y}):
(x,\mathbf{y})\in\overline{D_{S,2\delta}}\right\}
$$
and
$$
S_{\mathrm{side}}
=\left\{(x+\mathbf i t,\mathbf{y}):
(x,\mathbf{y})\in\partial D_{S,2\delta},\quad 0\le t\le\ell_n\right\},
$$
respectively.

\begin{figure}[!ht]
\centering
\includegraphics[width=.94\linewidth]{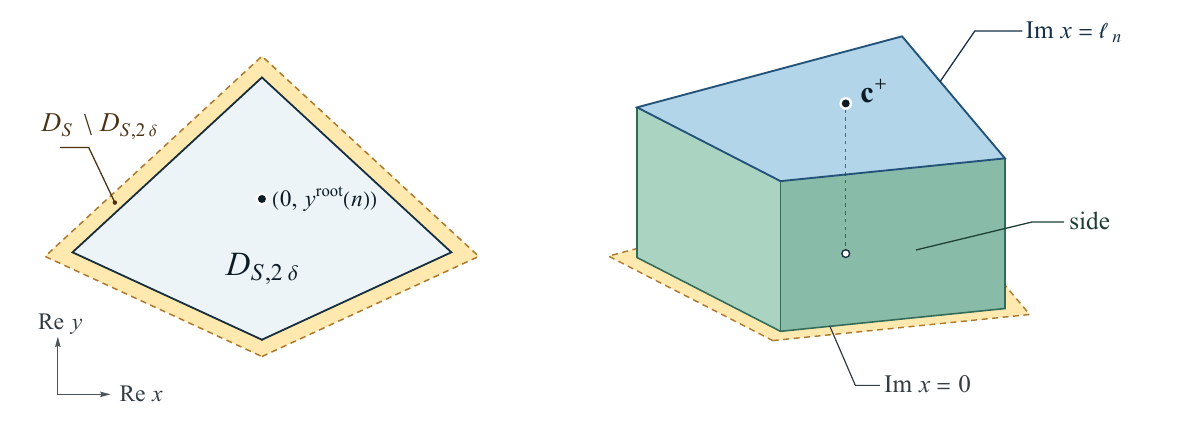}
\caption{The real domain $D_{S,2\delta}$ and its translation to
$\operatorname{Im}x=\ell_n$. }
\label{fig:one-edge-contour-shift}
\end{figure}

We next compare $V_r^+$ with $V^+$. By
Proposition~\ref{prop:double-sine-asymptotics}(1),
$$
\frac{2\pi\mathbf i}{r}
\log S_b\!\left(\frac{z}{\pi b}+\frac b2\right)
=\frac12L(z)+O\!\left(\frac1{r^2}\right).
$$
Replacing $z$ by $z-\frac{\pi}{r}$ and using
$L'(z)=2\mathbf i\log(2\sin z)$, we obtain
$$
\frac{2\pi\mathbf i}{r}\log S_b\!\left(\frac{z}{\pi b}\right)
=\frac12L(z)-\frac{\pi\mathbf i}{r}\log(2\sin z)
+O\!\left(\frac1{r^2}\right).
$$
Hence
\begin{equation}\label{eq:one-edge-quantum-phase-expansion}
\begin{aligned}
V_r^+(x,\mathbf{y})
={}&V^+(x,\mathbf{y})-\frac{\pi\mathbf i}{r}\sum\limits_{i=1}^n
\left[2\log(2\sin y_i)
+\log\!\bigl(2\sin(y_i+x)\bigr)
+\log\!\bigl(2\sin(y_i-x)\bigr)\right]\\
&+\frac{2\pi x}{r}+O\!\left(\frac1{r^2}\right).
\end{aligned}
\end{equation}
The estimate is uniform on a neighborhood of the deformation, since
the real parts of the arguments of $L$ lie in a compact subinterval
of $(0,\pi)$. Writing the remainder in
\eqref{eq:one-edge-quantum-phase-expansion} as
$\frac{2\pi\mathbf i}{r^2}\upsilon_r(x,\mathbf{y})$, we have
$$
2(1-e^{2\mathbf i x})
e^{\frac{r}{2\pi\mathbf i}V_r^+(x,\mathbf{y})}
=g(x,\mathbf{y})e^{r f_r(x,\mathbf{y})},
$$
where
$$
f_r(x,\mathbf{y})
=\frac{V^+(x,\mathbf{y})}{2\pi\mathbf i}
+\frac{\upsilon_r(x,\mathbf{y})}{r^2},
\qquad
g(x,\mathbf{y})
=-4\mathbf i\sin x
\prod\limits_{i=1}^{n}
\frac{1}{4\sin y_i\sqrt{\sin(y_i+x)\sin(y_i-x)}}.
$$
The square roots are positive on the real domain and are continued
along the deformation.

On $S_{\mathrm{top}}$, Lemma~\ref{lem:one-edge-shifted-potential}(1) gives
$$
\operatorname{Im}V^+(x+\mathbf i\ell_n,\mathbf{y})
\le\operatorname{Im}V^+(\mathbf i\ell_n,\mathbf{y}),
$$
with equality only at $x=0$. By part (2) of the same lemma, each term
$\operatorname{Im}\tau(\mathbf i\ell_n,y_i)$ attains its unique maximum
at $y_i=y^{\mathrm{root}}(n)$. Thus the unique maximum on $S_{\mathrm{top}}$ is attained
at $\mathbf{c}^+$, where its value is $\operatorname{Vol}(M_n)$ by
Proposition~\ref{prop:geometric-critical-points}.

On $S_{\mathrm{side}}$, first suppose that $|x|\le\frac{\pi}{4}$.
By Lemma~\ref{lem:one-edge-shifted-potential}(3),
$$
t\longmapsto\operatorname{Im}V^+(x+\mathbf i t,\mathbf{y})
=-\pi t+\sum\limits_{i=1}^{n}\operatorname{Im}\tau(x+\mathbf i t,y_i)
$$
is convex. It suffices to consider $t=0$ and $t=\ell_n$.
At $t=0$, the proof of Lemma~\ref{lem:one-edge-poisson} gives
$$
\operatorname{Im}V^+(x,\mathbf{y})
=\sum\limits_{i=1}^{n}\operatorname{Im}\tau(x,y_i)
\le\operatorname{Vol}(M_n)-\varepsilon',
\qquad (x,\mathbf{y})\in\partial D_{S,2\delta}.
$$
At $t=\ell_n$, the same conclusion follows from the uniqueness of the
maximum on $S_{\mathrm{top}}$, since $\mathbf{c}^+$ lies in its interior.

If $|x|\ge\frac{\pi}{4}$, Lemma~\ref{lem:one-edge-shifted-potential}(1) gives
$$
\operatorname{Im}V^+(x+\mathbf i t,\mathbf{y})
\le\operatorname{Im}V^+\!\left(\frac{\pi}{4}+\mathbf i t,\mathbf{y}\right).
$$
Again, part (3) reduces the estimate to the two endpoints. At $t=\ell_n$,
the real part of the first coordinate is $\frac{\pi}{4}$, so the point
differs from $\mathbf{c}^+$. At $t=0$, each $y_i$ lies in
$[\frac{5\pi}{8},\frac{3\pi}{4}]$. By \eqref{eq:one-edge-classical-im}
and the concavity of $\Lambda$ on $[0,\frac{\pi}{2}]$,
$$
\operatorname{Im}\tau\!\left(\frac{\pi}{4},y_i\right)
\le4\Lambda\!\left(\frac{5\pi}{16}\right)
+3\Lambda\!\left(\frac{\pi}{6}\right)
<3.1<\frac{\operatorname{Vol}(M_2)}{2}
\le\frac{\operatorname{Vol}(M_n)}{n}.
$$
It follows by compactness that
$$
\sup\limits_{S_{\mathrm{side}}}\operatorname{Im}V^+
<\operatorname{Vol}(M_n).
$$

By Cauchy's theorem, we may deform $D_{S,2\delta}$ to the oriented
surface $S_{\mathrm{top}}\cup S_{\mathrm{side}}$. The contribution
from $S_{\mathrm{side}}$ is exponentially smaller, and
$\operatorname{Im}V^+$ attains its unique maximum on $S_{\mathrm{top}}$
at $\mathbf{c}^+$. Set
$$
f=\frac{V^+}{2\pi i}.
$$
Choose $\rho>0$ sufficiently small and set
$$
\begin{aligned}
D&=\left\{
c^++u+iv:
u,v\in\mathbb R^{n+1},\
\|u\|<2\rho,\ \|v\|<\rho
\right\},\\[4pt]
S&=\left\{
c^++u:
u\in\mathbb R^{n+1},\
\|u\|\le\rho
\right\}
\subset S_{\mathrm{top}}.
\end{aligned}
$$
Since $\operatorname{Re}f$ attains its unique maximum on
$S_{\mathrm{top}}$ at $c^+$, the contribution from
$S_{\mathrm{top}}\setminus S$ is exponentially smaller.
We therefore apply Proposition~5.4 to $S$ with this choice of $D$.

By Lemma~5.3, the real Hessian of $\operatorname{Re}f$
in the horizontal directions at $c^+$ is negative definite.
Since $f$ is holomorphic, after decreasing $\rho$ if necessary,
$\operatorname{Re}f$ is strictly concave in the horizontal
directions and strictly convex in the imaginary directions on $D$.
The homotopy
$$
c^++u+iv\longmapsto c^++u+i(1-t)v,
\qquad 0\le t\le1,
$$
followed by radial retraction onto $S$, gives a deformation
retraction of
$$
\left\{
z\in D:
\operatorname{Re}f(z)<\operatorname{Re}f(c^+)
\right\}
$$
onto $S\setminus\{c^+\}$, as required in condition~(iii).

Using $4\cosh\ell_n\sin^2y^{\mathrm{root}}(n)=1+\cosh\ell_n$ and
Proposition~\ref{prop:hessian-torsion}, we obtain
$$
\frac{g(\mathbf{c}^+)}
{\sqrt{\det\!\left(\mathbf i\,\operatorname{Hess}V^+(\mathbf{c}^+)\right)}}
=\frac{2}{\sqrt{\left|\operatorname{Tor}(DM_n;\operatorname{Ad}\rho)\right|}}.
$$
Here the Gaussian square root is positive.
Finally, Proposition~\ref{prop:saddle-point-method} in $n+1$ variables
gives the  factor
$\left(\frac{4\pi^2}{r}\right)^{\frac{n+1}{2}}$,
and Proposition~\ref{prop:geometric-critical-points} gives
$V^+(\mathbf{c}^+)=\mathbf i\,\operatorname{Vol}(M_n)$. Substituting these into
the saddle point formula proves \eqref{eq:one-edge-leading-asymptotic}.
\end{proof}

\subsection{Estimates for the remaining Fourier coefficients}

We now estimate the remaining Fourier coefficients in \eqref{eq:one-edge-poisson}.
We group $f_r^+(p,\mathbf q)$ with $f_r^-(p-1,\mathbf q)$, so the four
leading terms correspond to $p=0,1$ and $\mathbf q=\mathbf 0$.

\begin{proposition}\label{prop:one-edge-remaining-modes}
For fixed $n\ge2$, there exists $\varepsilon>0$ such that, as
$r\to\infty$ through odd integers,
\begin{equation}\label{eq:one-edge-remaining-estimate}
\sum\limits_{\substack{(p,\mathbf{q})\in\mathbb Z^{n+1}\\
\mathbf{q}\ne\mathbf{0}\ \mathrm{or}\ p\notin\{0,1\}}}
\left(|f_r^+(p,\mathbf{q})|+|f_r^-(p-1,\mathbf{q})|\right)
=O\!\left(e^{\frac{r}{2\pi}(\operatorname{Vol}(M_n)-\varepsilon)}\right).
\end{equation}
\end{proposition}

\begin{proof}
We first integrate by parts on the real domain. For a smooth function $f$, write
$$
\Delta f=f''_{xx}+\sum\limits_{i=1}^{n}f''_{y_i y_i}.
$$
Since
$$
\left(1-\frac{\Delta}{r^2}\right)
e^{-\mathbf i r(px+\mathbf{q}\cdot\mathbf{y})}
=(1+p^2+\|\mathbf{q}\|^2)
e^{-\mathbf i r(px+\mathbf{q}\cdot\mathbf{y})}
$$
and $\psi$ has compact support in $D_S$, integrating by parts $n+1$
times gives
$$
f_r^\pm(p,\mathbf{q})=\frac{(-1)^p}{(1+p^2+\|\mathbf{q}\|^2)^{n+1}}\int_{D_S}\phi_r(x,\mathbf{y})\,e^{\frac{r}{2\pi\mathbf i}\left(V_r^\pm(x,\mathbf{y})+2\pi p x+2\pi\mathbf{q}\cdot\mathbf{y}\right)}\,dx\,d\mathbf{y},
$$
where, for each choice of sign, $\phi_r$ is determined by
$$
\left(1-\frac{\Delta}{r^2}\right)^{n+1}
\left(\psi\,e^{\frac{r}{2\pi\mathbf i}V_r^\pm}\right)
=\phi_r\,e^{\frac{r}{2\pi\mathbf i}V_r^\pm}.
$$
By \eqref{eq:one-edge-quantum-phase-expansion},
\eqref{eq:quantum-potential-symmetry}, and Cauchy's estimates, $|\phi_r|\le C$ on
$\operatorname{supp}\psi$, uniformly in $r$. Since $\psi=1$ on
$D_{S,2\delta}$, the restriction of $\phi_r$ to $D_{S,2\delta}$
extends holomorphically to $D_S+\mathbf i\mathbb R^{n+1}$ and is
uniformly bounded in $r$ on compact subsets of this tube. In particular,
$$
\sup\limits_{\substack{(x,\mathbf{y})\in\overline{D_{S,2\delta}}\\
0\le t\le\ell_n}}
|\phi_r(x+\mathbf i t,\mathbf{y})|\le C.
$$

Since $V^-(x,\mathbf y)=V^+(x,\mathbf y)+2\pi x$, the linear term in $x$
is $\pi(2p-1)x$ for both $f_r^+(p,\mathbf q)$ and $f_r^-(p-1,\mathbf q)$
in the classical limit. The identity $f_r^+(p,\mathbf q)=f_r^-(-p,\mathbf q)$
interchanges the two terms under $p\mapsto1-p$, which pairs $p\le0$ with $p\ge1$.
Thus
\begin{equation}\label{eq:one-edge-remaining-pairs}
\begin{aligned}
&\sum\limits_{\substack{(p,\mathbf{q})\in\mathbb Z^{n+1}\\
\mathbf{q}\ne\mathbf{0}\ \mathrm{or}\ p\notin\{0,1\}}}
\left(|f_r^+(p,\mathbf{q})|+|f_r^-(p-1,\mathbf{q})|\right)\\
=&2
\sum\limits_{\substack{(p,\mathbf{q})\in\mathbb Z^{n+1}\\p\le-1}}
\left(|f_r^+(p,\mathbf{q})|+|f_r^-(p-1,\mathbf{q})|\right)
+2\sum\limits_{\mathbf{q}\in\mathbb Z^n\setminus\{\mathbf{0}\}}
\left(|f_r^+(0,\mathbf{q})|+|f_r^-(-1,\mathbf{q})|\right).
\end{aligned}
\end{equation}
On the right-hand side, $p\le0$. Thus for real
$x,\mathbf{y}$ and $0\le t\le\ell_n$,
\begin{equation}\label{eq:one-edge-upper-phase}
\begin{aligned}
&\operatorname{Im}\!\left(V^+(x+\mathbf i t,\mathbf{y})
+2\pi p(x+\mathbf i t)+2\pi\mathbf{q}\cdot\mathbf{y}\right)=\operatorname{Im}V^+(x+\mathbf i t,\mathbf{y})
+2\pi p t
\le\operatorname{Im}V^+(x+\mathbf i t,\mathbf{y}).
\end{aligned}
\end{equation}

On $D_S\setminus D_{S,2\delta}$, the Fourier factor has absolute value
$1$, so the real collar estimate in the proof of
Proposition~\ref{prop:one-edge-leading-terms} applies uniformly in
$p,\mathbf{q}$. By Cauchy's theorem, we deform $D_{S,2\delta}$ to
$S_{\mathrm{top}}\cup S_{\mathrm{side}}$ by translating $x$ to
$x+\mathbf i\ell_n$. On $S_{\mathrm{side}}$,
\eqref{eq:one-edge-upper-phase} and the same proposition give a uniform
gap below $\operatorname{Vol}(M_n)$. It remains to estimate the
integrals on $S_{\mathrm{top}}$.

\medskip
\noindent\textit{Case 1.} Suppose $p\le-1$. Taking $t=\ell_n$
in \eqref{eq:one-edge-upper-phase}, we obtain
$$
\sup\limits_{(x,\mathbf{y})\in S_{\mathrm{top}}}
\operatorname{Im}\!\left(V^+(x,\mathbf{y})+2\pi p x
+2\pi\mathbf{q}\cdot\mathbf{y}\right)
\le\operatorname{Vol}(M_n)-2\pi\ell_n.
$$
This gives the required gap for the first sum in
\eqref{eq:one-edge-remaining-pairs}.

\medskip
\noindent\textit{Case 2.} Suppose $p=0$ and
$\mathbf{q}\ne\mathbf{0}$. For both terms in the second sum in
\eqref{eq:one-edge-remaining-pairs}, the exponential factor in the integrand
has growth rate
$$
\operatorname{Im}\!\left(V^+(x,\mathbf{y})+2\pi\mathbf{q}\cdot\mathbf{y}\right).
$$
Choose a small closed ball $B$ centered at $\mathbf{c}^+$ in the interior of
$S_{\mathrm{top}}$. Since $(V^+)'_{y_i}(\mathbf{c}^+)=0$ for every $i$,
we may choose $B$ so that
$$
\left|\operatorname{Re}(V^+)'_{y_i}\right|\le\pi,
\qquad 1\le i\le n,
$$
on a complex neighborhood of $B$. By the uniqueness of the maximum
on $S_{\mathrm{top}}$,
$$
\sup\limits_{S_{\mathrm{top}}\setminus\operatorname{int}B}
\operatorname{Im}V^+\le\operatorname{Vol}(M_n)-\varepsilon'
$$
for some $\varepsilon'>0$.

Take $i$ such that $q_i\ne0$. Keeping $\operatorname{Im}x=\ell_n$,
translate $B$ by
$$
y_i\longmapsto y_i-\mathbf i u\,\operatorname{sgn}(q_i),
\qquad 0\le u\le h,
$$
where $h>0$ is sufficiently small and independent of $\mathbf{q}$.
Along this translation,
$$
\left[\operatorname{Im}\bigl(V^++2\pi\mathbf{q}\cdot\mathbf{y}\bigr)\right]'_u
=-\operatorname{sgn}(q_i)
\operatorname{Re}(V^+)'_{y_i}
-2\pi|q_i|\le-\pi.
$$
Hence this growth rate is at most
$\operatorname{Vol}(M_n)-\pi h$ on the translated ball and at most
$\operatorname{Vol}(M_n)-\varepsilon'$ on the sides swept out by
$\partial B$. By Cauchy's theorem, this gives the required gap for the
second sum in \eqref{eq:one-edge-remaining-pairs}.

Finally, \eqref{eq:one-edge-quantum-phase-expansion} and
\eqref{eq:quantum-potential-symmetry} give $V_r^\pm=V^\pm+O(\frac1r)$ uniformly on the
compact deformation regions above. Since $\phi_r$ is uniformly
bounded there, the two cases give, for the indices on the right-hand side
of \eqref{eq:one-edge-remaining-pairs},
$$
|f_r^+(p,\mathbf{q})|+|f_r^-(p-1,\mathbf{q})|
\le\frac{C}{(1+p^2+\|\mathbf{q}\|^2)^{n+1}}
e^{\frac{r}{2\pi}(\operatorname{Vol}(M_n)-\varepsilon)}
$$
for some $\varepsilon>0$ independent of $p,\mathbf{q},r$.
Summing over \eqref{eq:one-edge-remaining-pairs} and using
$$
\sum\limits_{(p,\mathbf{q})\in\mathbb Z^{n+1}}
\frac{1}{(1+p^2+\|\mathbf{q}\|^2)^{n+1}}<\infty
$$
proves \eqref{eq:one-edge-remaining-estimate}.
\end{proof}
\subsection{Proof of the main theorem}

By \eqref{eq:one-edge-poisson}, we have
$$
\mathrm{TV}_r(M_n)
=
\frac{(-1)^{\frac{(r+1)(n+1)}{2}}\,r}{8\pi^{n+1}}
\sin^{n-1}\frac{2\pi}{r}
\sum\limits_{\sigma\in\{+,-\}}
\sum\limits_{(p,\mathbf{q})\in\mathbb Z^{n+1}}
f_r^\sigma(p,\mathbf{q})
+
O\!\left(e^{\frac{r}{2\pi}(\operatorname{Vol}(M_n)-\varepsilon)}\right).
$$

By Proposition~\ref{prop:one-edge-leading-terms}, the four leading
terms satisfy
$$
f_r^+(0,\mathbf{0})+f_r^+(1,\mathbf{0})
+f_r^-(0,\mathbf{0})+f_r^-(-1,\mathbf{0})
=
2\left(\frac{4\pi^2}{r}\right)^{\frac{n+1}{2}}
\frac{e^{\frac{r}{2\pi}\operatorname{Vol}(M_n)}}
{\sqrt{|\operatorname{Tor}(DM_n)|}}
\left(1+O\!\left(\frac1r\right)\right).
$$
By Proposition~\ref{prop:one-edge-remaining-modes}, the remaining
terms satisfy
$$
\sum\limits_{\substack{(p,\mathbf{q})\in\mathbb Z^{n+1}\\
\mathbf{q}\ne\mathbf{0}\ \mathrm{or}\ p\notin\{0,1\}}}
\left(|f_r^+(p,\mathbf{q})|+|f_r^-(p-1,\mathbf{q})|\right)
=
O\!\left(e^{\frac{r}{2\pi}(\operatorname{Vol}(M_n)-\varepsilon)}\right).
$$

Substituting these estimates into \eqref{eq:one-edge-poisson} and using
$\sin\frac{2\pi}{r}
=\frac{2\pi}{r}\left(1+O\!\left(\frac{1}{r^2}\right)\right)$,
we obtain
$$
\mathrm{TV}_r(M_n)
=
(-1)^{\frac{(r+1)(n+1)}{2}}
\frac{(4\pi)^{n-1}}
{r^{\frac{3(n-1)}{2}}\sqrt{|\operatorname{Tor}(DM_n)|}}
\,e^{\frac{r}{2\pi}\operatorname{Vol}(M_n)}
\left(1+O\!\left(\frac1r\right)\right).
$$
The volume limit follows by taking absolute values and logarithms.
\hfill$\square$

\appendix
\section{Proofs}\label{app:sixj-proofs}

\subsection{Special function identities}\label{app:special-function-identities}

We first prove \eqref{eq:factorial-double-sine}. For periods $b,b^{-1}$, the conventions of \cite[(19)]{Narukawa2004}
and \cite[(2.1)]{MYComplex} give
\[
 S_2(w\mid b,b^{-1})=\frac{1}{S_b(w)}.
\]
For $\operatorname{Re}b>0$ and $\operatorname{Im}(b^2)>0$, Narukawa's
product formula~\cite[Corollary~6]{Narukawa2004} reads
\begin{equation}\label{eq:narukawa-product}
 S_2(w\mid b,b^{-1})=
 \exp\!\left(\frac{\pi\mathbf i}{2}B_{22}(w\mid b,b^{-1})\right)
 \frac{\bigl(e^{2\pi\mathbf i bw};q^2\bigr)_\infty}
 {\bigl(e^{\frac{2\pi\mathbf i}{b}(w-b^{-1})};e^{-\frac{2\pi\mathbf i}{b^2}}\bigr)_\infty},
\end{equation}
where $B_{22}(w\mid b,b^{-1})=w^2-Qw+\frac{b^2+b^{-2}+3}{6}$.
When $w$ changes from $b$ to $(z+1)b$, the first Pochhammer argument
changes from $q^2$ to $q^{2(z+1)}$ and the second is multiplied by
$e^{2\pi\mathbf i z}$. The polynomial difference
$B_{22}(b\mid b,b^{-1})-B_{22}((z+1)b\mid b,b^{-1})$ is
$z-b^2z(z+1)$. Taking the quotient of
\eqref{eq:narukawa-product} at those two arguments and using the
complex-order Pochhammer symbol of Section~\ref{subsec:quantum-sixj} gives
\[
 \begin{aligned}
 \frac{S_b((z+1)b)}{S_b(b)}
 &=\frac{S_2(b\mid b,b^{-1})}{S_2((z+1)b\mid b,b^{-1})}\\
 &=\exp\!\left(\frac{\pi\mathbf i}{2}\bigl(z-b^2z(z+1)\bigr)\right)
   \frac{(q^2;q^2)_z}
   {\bigl(e^{-\frac{2\pi\mathbf i}{b^2}};e^{-\frac{2\pi\mathbf i}{b^2}}\bigr)_{-b^2z}} .
 \end{aligned}
\]
Multiply by $(2\sin(\pi b^2))^{-z}$. Since
$2\sin(\pi b^2)=\mathbf i q^{-1}(1-q^2)$, the right-hand side becomes
\eqref{eq:pochhammer-factorial}, which proves \eqref{eq:factorial-double-sine}.
For $\operatorname{Im}(b^2)>0$ the logarithm used in the power on the
right is $\frac{\pi\mathbf i}{2}-\pi\mathbf i b^2+\operatorname{Log}(1-q^2)$.
It has a boundary value at every real $b>0$ with $b^2\notin\mathbb Z$.
With this branch \eqref{eq:factorial-double-sine} gives the continuation
there as a meromorphic function of $z$.

For the first identity in \eqref{eq:factorial-complement}, note that
$\frac{Q_4}{2}=b^{-2}-1$, so the double sine arguments $(z+1)b$ and
$\bigl(\frac{Q_4}{2}-z+1\bigr)b$ sum to $Q=b+b^{-1}$.
By \eqref{eq:factorial-double-sine} and $S_b(w)S_b(Q-w)=1$, the product
$[z]!\left[\frac{Q_4}{2}-z\right]!$ is constant in $z$. At $z=0$ it equals
$\left[\frac{Q_4}{2}\right]!$.
For the second identity, the shift equation $S_b(w+b)=2\sin(\pi bw)S_b(w)$
gives $[z+1]!=[z+1][z]!$. Iterating it $n$ times gives
\[
 \frac{[Q_4+1]!}{[Q_4+1-n]!}
 =\prod\limits_{j=1}^n[2b^{-2}-j]=(-1)^n[n]!,
\]
since $q=e^{\pi\mathbf i b^2}$ implies $[2b^{-2}-j]=-[j]$.\qed

\subsection{Square root conventions and Racah normalization}\label{app:root-signs}

Recall from \eqref{eq:factorial-kernel} that
\[
 U_b(\mathbf a;k)=
 -\frac{1}{[Q_4+1]!}\,
 \frac{\sqrt{\displaystyle\prod\limits_{i,j=1}^4[Q_j-T_i]!}}
 {\displaystyle\prod\limits_{i=1}^4[k-T_i]!\prod\limits_{j=1}^4[Q_j-k]!}.
\]

For an admissible triple $(a_1,a_2,a_3)$, define the triangle factor
\begin{equation}\label{eq:triangle-factor}
\Delta(a_1,a_2,a_3)=
\sqrt{\frac{
\left[\frac{a_1+a_2-a_3}{2}\right]!\,
\left[\frac{a_2+a_3-a_1}{2}\right]!\,
\left[\frac{a_3+a_1-a_2}{2}\right]!}
{\left[\frac{a_1+a_2+a_3}{2}+1\right]!}}.
\end{equation}
Take the principal square root of the quotient, with
$\sqrt{x}=\mathbf i\sqrt{-x}$ for $x<0$, following
\cite[Section~3.4]{YangRelative}. For an admissible $6$-tuple
$\mathbf a$, the Racah prefactor is
\[
C(\mathbf{a})=\mathbf i^{-\sum\limits_{i=1}^6 a_i}
\Delta(a_1,a_2,a_3)\Delta(a_1,a_5,a_6)
\Delta(a_2,a_4,a_6)\Delta(a_3,a_4,a_5).
\]
The numerator in \eqref{eq:factorial-kernel}
is a single square root. For the original tuple $\mathbf a$ we take
$([Q_4+1]!)^2C(\mathbf a)$. For both reflected sums in
Theorem~\ref{thm:hyperideal-finite-sums} we take the root
\eqref{eq:reflected-root} with the indices relabeled for the chosen reflection and the original $C(\mathbf a)$.

For an admissible $\mathbf a$ the 12
differences $Q_j-T_i$ with $j\le3$ are the triangle differences in
\eqref{eq:triangle-factor}. Applying \eqref{eq:factorial-complement}
with $n=T_i+1$ shows that the radicand in \eqref{eq:factorial-kernel}
is $([Q_4+1]!)^4 C(\mathbf a)^2$, since $\sum\limits_{i=1}^4 T_i=\sum\limits_{i=1}^6 a_i$.
A second application, with $n=k+1$, gives the usual Racah summand
for every integer $k$ in the bounds of \eqref{eq:quantum-sixj}:
\begin{equation}\label{eq:racah-term}
 U_b(\mathbf a;k)=C(\mathbf a)
 \frac{(-1)^k[k+1]!}
 {\displaystyle\prod\limits_{i=1}^4[k-T_i]!\prod\limits_{j=1}^3[Q_j-k]!}.
\end{equation}
Thus Definition~\ref{def:quantum-sixj} agrees with the standard
normalization.

At $b^2=\frac{2}{r}$, an $r$-hyperideal tuple satisfies
$0\le Q_j-T_i<\frac{r-2}{2}$ for $1\le i,j\le4$. Thus
$(Q_j-T_i+1)b\in(0,b+b^{-1})$, where the individual square root
$\sqrt{S_b}$ is positive. The single numerator root in \eqref{eq:root-kernel}
has the sign of $C(\mathbf a)$, which is preserved under the transport
in Appendix~\ref{app:reflected-sums}.

\subsection{Proof of Proposition~\ref*{prop:double-sine-sixj}}\label{app:double-sine-proof}

For $b^{-2}\notin\mathbb Z$, we prove the double sine formula \eqref{eq:root-kernel},
\[
 U_b(\mathbf a;k)=
 -\frac{b^2\sin(\pi b^2)}{\sin\!\left(\frac{\pi}{b^2}\right)}
 \frac{\sqrt{\displaystyle\prod\limits_{i,j=1}^4 S_b((Q_j-T_i+1)b)}}
 {\displaystyle\prod\limits_{i=1}^4 S_b((k-T_i+1)b)
  \prod\limits_{j=1}^4 S_b((Q_j-k+1)b)} .
\]

Substitute \eqref{eq:factorial-double-sine} into \eqref{eq:factorial-kernel}. The factors $S_b(b)$ cancel, and
$\sum\limits_{j=1}^4Q_j-\sum\limits_{i=1}^4T_i=Q_4$ leaves the factor
$(2\sin(\pi b^2))^{-Q_4}$ in the quotient. The
\hyperref[sb:functional-identities]{identities in Section~\ref*{subsec:special-functions}}
give
\[
 S_b\!\left(\frac{2}{b}\right)
 =2\sin\!\left(\frac{\pi}{b^2}\right)S_b\!\left(\frac{1}{b}\right)
 =\frac{2\sin\!\left(\frac{\pi}{b^2}\right)}{b}.
\]
Substituting $z=Q_4+1$ in \eqref{eq:factorial-double-sine} now gives
the coefficient $-\frac{b^2\sin(\pi b^2)}{\sin\!\left(\frac{\pi}{b^2}\right)}$ in \eqref{eq:root-kernel}.\qed

\subsection{Proof of Theorem~\ref*{thm:hyperideal-finite-sums}}\label{app:reflected-sums}

We prove \eqref{eq:reflected-integer} and \eqref{eq:reflected-shifted} for
$F_4(\mathbf a)=(a_1,a_2,a_3,Q_4-a_4,a_5,a_6)$. This reflection fixes
$T_1,T_2,Q_3$ and adds $\delta=\frac{Q_4}{2}-a_4$ to $T_3,T_4,Q_1,Q_2$.
Thus the two identities to be proved in this case are
\[
 \begin{aligned}
 \left|\begin{matrix}a_1&a_2&a_3\\a_4&a_5&a_6\end{matrix}\right|
 &=\sum\limits_{k=\max\{T_1,T_2\}}^{Q_3} U_b(F_4(\mathbf a);k),\\
 &=\sum\limits_{k'=\max\{T_3,T_4\}}^{\min\{Q_1,Q_2\}} U_b(F_4(\mathbf a);k'+\delta),
 \end{aligned}
 \qquad k,k'\in\mathbb Z .
\]
The other $F_i$ follow by relabeling.

We use the single-sum identities for generic $q$ from
\cite{SmirnovTolstoyKharitonov1991,AsherovaEtAlCzech1996,AsherovaEtAlPAN1996},
as recorded by Ali{\v{s}}auskas~\cite{Alisauskas2000}. For the choice $T_1,T_2,Q_3$, put
$t=T_3+T_4=Q_3+a_4$. The common prefactor in our normalization is
\[
 P=C(\mathbf a)\,
 \frac{[T_3+1]![T_4+1]!}
 {[Q_3-T_3]![Q_3-T_4]!\displaystyle\prod\limits_{i,j=1}^2[Q_j-T_i]!}.
\]
The first formula, \cite[(2.2a)]{Alisauskas2000}, reads
\begin{equation}\label{eq:al-first-sum}
 \left|\begin{matrix}a_1&a_2&a_3\\a_4&a_5&a_6\end{matrix}\right|
 =P\!\sum\limits_{k=\max\{T_1,T_2\}}^{Q_3}
 \frac{(-1)^k[k+1]!\displaystyle\prod\limits_{j=1}^2[k+a_4-Q_j]!}
 {[Q_3-k]!\displaystyle\prod\limits_{i=1}^2[k-T_i]!
              \prod\limits_{i=3}^4[k+1+a_4-T_i]!}.
\end{equation}
The second formula, \cite[(2.1a)]{Alisauskas2000}, reads
\begin{equation}\label{eq:al-second-sum}
 \left|\begin{matrix}a_1&a_2&a_3\\a_4&a_5&a_6\end{matrix}\right|
 =P\!\sum\limits_{k'=\max\{T_3,T_4\}}^{\min\{Q_1,Q_2\}}
 \frac{(-1)^{k'}[k'-a_4]!\displaystyle\prod\limits_{i=1}^2[a_4+T_i-k']!}
 {[t+1-k']!\displaystyle\prod\limits_{i=3}^4[k'-T_i]!
               \prod\limits_{j=1}^2[Q_j-k']!}.
\end{equation}

Substituting $F_4(\mathbf a)$ in the definition \eqref{eq:factorial-kernel} gives
\begin{equation}\label{eq:f4-kernel}
 U_b(F_4(\mathbf a);k)=-
 \frac{\sqrt{\displaystyle
   \prod\limits_{\substack{i=1,2\\j=3,4}}[Q_j-T_i]!\,
   \prod\limits_{\substack{i=3,4\\j=1,2}}[Q_j-T_i]!\,
   \prod\limits_{i,j=1}^2[Q_j-T_i+\delta]!\,
   \prod\limits_{i,j=3}^4[Q_j-T_i-\delta]!}}
 {[Q_4+1]![Q_3-k]![Q_4-k]!\displaystyle
   \prod\limits_{i=1}^2[k-T_i]!\,
   \prod\limits_{i=3}^4[k-T_i-\delta]!\,
   \prod\limits_{j=1}^2[Q_j+\delta-k]!}.
\end{equation}
Write $R_4$ for the numerator square root. In its last two products
the arguments are $\frac{Q_4}{2}-d$, up to permutation, where
$d=Q_j-T_i$ for $i,j\in\{1,2\}$ or $i,j\in\{3,4\}$. This follows from
$Q_1+Q_2=T_1+T_2+a_4$, $Q_3=T_3+T_4-a_4$ and the definition of $\delta$.
Since $\left[\frac{Q_4}{2}-d\right]!=\frac{\left[\frac{Q_4}{2}\right]!}{[d]!}$, each replacement multiplies the
radicand by $\frac{\left[\frac{Q_4}{2}\right]!}{([d]!)^2}$. Transporting the original root and
using \eqref{eq:factorial-complement} at $n=T_3+1,T_4+1$ gives
\begin{equation}\label{eq:reflected-root}
 R_4=([Q_4+1]!)^2C(\mathbf a)\frac{\left(\left[\frac{Q_4}{2}\right]!\right)^4}{\prod\limits_d[d]!}
    =(-1)^t\left(\left[\frac{Q_4}{2}\right]!\right)^4P .
\end{equation}
For an $r$-hyperideal tuple, $0\le d<\frac{Q_4}{2}$ for every such $d$, so the
factor multiplying the original root is positive.

For integer $k$, applying \eqref{eq:factorial-complement} to the four
denominator factors involving $\delta$ in \eqref{eq:f4-kernel} gives
\[
 \begin{aligned}
 \frac{1}{[k-T_i-\delta]!}
 &=\frac{(-1)^{k+1+a_4-T_i}[Q_4+1]!}{\left[\frac{Q_4}{2}\right]![k+1+a_4-T_i]!}
 &&(i=3,4),\\
 \frac{1}{[Q_j+\delta-k]!}
 &=\frac{[k+a_4-Q_j]!}{\left[\frac{Q_4}{2}\right]!} &&(j=1,2).
 \end{aligned}
\]
The two signs in the first row multiply to $(-1)^t$. For the factor
$[Q_4-k]!$ in the denominator, use
$\frac{1}{[Q_4-k]!}=(-1)^{k+1}\frac{[k+1]!}{[Q_4+1]!}$.
Together with \eqref{eq:reflected-root}, these substitutions cancel
$\left(\left[\frac{Q_4}{2}\right]!\right)^4$ and the factors $[Q_4+1]!$. The remaining prefactor is
$P$, and the sign is $(-1)^k$, giving the summand of \eqref{eq:al-first-sum}.
Summing from $\max(T_1,T_2)$ to $Q_3$ proves \eqref{eq:reflected-integer}.

For the second sum set $k=k'+\delta$, with $k'\in\mathbb Z$.
The four denominator factors involving $\delta$ in \eqref{eq:f4-kernel}
now have  arguments $k'-T_3,k'-T_4,Q_1-k',Q_2-k'$. Applying
\eqref{eq:factorial-complement} to the other $k$-dependent factors gives
\[
 \begin{aligned}
 \frac{1}{[k-T_i]!}&=\frac{[a_4+T_i-k']!}{\left[\frac{Q_4}{2}\right]!} &&(i=1,2),\\
 \frac{1}{[Q_4-k]!}&=\frac{[k'-a_4]!}{\left[\frac{Q_4}{2}\right]!},\\
 \frac{1}{[Q_3-k]!}&=\frac{(-1)^{t+1-k'}[Q_4+1]!}{\left[\frac{Q_4}{2}\right]![t+1-k']!}.
 \end{aligned}
\]
The sign from $[Q_3-k]!$, the sign $(-1)^t$ in $R_4$ and the minus sign
in \eqref{eq:f4-kernel} give $(-1)^{k'}$. Cancelling the same constants
gives the summand of \eqref{eq:al-second-sum}. Under $k=k'+\delta$ its
bounds become $\max(T_3,T_4)+\delta\le k\le\min(Q_1,Q_2)+\delta$.
Summing over this range proves \eqref{eq:reflected-shifted}.

The $G_{i'}$-substitutions follow since $G_{i'}$ permutes the $T$'s and
$Q$'s separately in \eqref{eq:racah-term} with the original prefactor $C(\mathbf a)$.\qed

\bibliographystyle{abbrv}
\bibliography{bibtex}

\end{document}